\documentclass{article}
\usepackage{latexsym,amssymb,amsmath,amsthm}
\usepackage[english]{babel}
\usepackage[twoside,dvips]{geometry}
\usepackage[all]{xy}
\newcommand{\R}{{\mathbb{R}}}
\newcommand{\C}{{\mathbb{C}}}

\newcommand{\n}{\noindent}
\newcommand{\mioplus}[1][]{\,\,\oplus_{#1}\,\,}

\def\H{\text{$\mathbb{H}$}}
\newcommand\SO{{\rm SO}}

\newcommand\SU{{\rm SU}}
\newcommand\Su{{\rm S}}
\def\Sp{{\rm Sp}}
\newcommand\U{{\rm U}}
\newcommand\SUr[2]{{{\rm S}({\rm U}({#1}) \times {\rm U}({#2}))}}
\newcommand\GL{{\rm GL}}
\newcommand\SL{{\rm SL}}

\newcommand\Irr{{\rm Irr}}

\newcommand\mf[1]{\mathfrak{#1}}
\newcommand\Spin{{\rm Spin }}

\newcommand\G{{\rm G}}

\newcommand\Ea{{\rm E}}

\newcommand\li[1]{\mathfrak{l}_{#1}}
\newcommand\h{\mathfrak{l}}
\newcommand\su{\mathfrak{su} }

\newcommand\asp{\mathfrak{sp}}

\newtheorem{thm}{Theorem}[section]
\newtheorem{prop}{Proposition}[section]
\newtheorem{defi}{Definition}[section]
\newtheorem{lem}{Lemma}[section]

\begin{document}
\pagestyle{plain}
\title{Coisotropic and Polar actions on Complex Grassmannians}
\author{Leonardo Biliotti and Anna Gori}
\date{}
\maketitle
%\thanks{Research partially supported by CNPq and FAPESP %(Brazil)}
%\email{gori@math.unifi.it}
%\address{Dipartimento di Matematica ``U.~Dini'', Universit\`a di
%Firenze, Viale Morgagni 67/a, I-50134, Florence, Ital
%\Tableofcontents
%\addcontentsline{toc}{section}{Introduction}
%\pagebreak
%\markright{Introduction}
%%%%%%%%%%%%%%%%%%%%%%%%%%%%%%%%%%%%%%%%%%%%%%%%%%%%%%%%%%%%%%%%%%%%%%%%%%%%%%%%%%%%%%%%%%%%%%%%%%%%%%%%%%%%%%%%%%%%%%%%%%%%%
%                                PARTE INTRODUTTIVA
%
%%%%%%%%%%%%%%%%%%%%%%%%%%%%%%%%%%%%%%%%%%%%%%%%%%%%%%%%%%%%%%%%%%%%%%%%%%%%%%%%%%%%%%%%%%%%%%%%%%%%%%%%%%%%%%%%%%%%%%%%%%%%
\begin{abstract}
\n
The main result of the paper is the complete classification of the compact connected Lie groups
acting coisotro\-pi\-cal\-ly
on complex Grassmannians. This is used to determine the polar actions
on the same manifolds.
\end{abstract}
\section{Introduction and Preliminaries}
The aim of this paper is to present the  classification of
 compact connected subgroups of $\SU(n)$ acting
coisotropically on the complex Grassmannians  $Gr(k,n)= \SU(n)/
\SUr{k}{n-k}$.\\ If $(M,g)$ is a compact K\"ahler manifold with
K\"ahler form $\omega$ and $K$ is a compact connected Lie subgroup
of its full isometry group, then the $K$-action is called {\em
coisotropic} or {\em multiplicity free} if the principal
$K$-orbits are coisotropic with respect to $\omega$ \cite{Hw}. Kac
\cite{Kac} and Benson and Ratcliff \cite{BR} have given  the
classification of linear multiplicity free
 representations, from which one has the full classification of coisotropic
actions on $Gr (k,n)$ for $k=1$, i.e. on the complex projective
space. It is therefore natural to investigate this kind of actions
on complex Grassmannians. Note that, in our analysis, we will
always assume $2 \leq k \leq \frac{n}{2}  ;$ if $k> \frac{n}{2}$ we refer to
the dual Grassmannian $Gr(k',n),$ where
$k'=n-k\leq \frac{n}{2}.$ Our main result is given in the following
\begin{thm} \label{risultato} Let $K$ be a compact
connected semisimple Lie subgroup of $\SU(n),$
acting non transitively on the complex Grassmannians
$M=Gr(k,n)= \SU(n)/ \SUr{k}{n-k}.$  Then $K$ acts coisotropically
on $M$ if and only if its Lie algebra $\mf{k}$ is  conjugate
to one of the Lie algebras  appearing in Table $1.$
\begin{small}
\begin{center}{\bf Table 1}
\end{center}
$$
\begin{array}{|c|c|c|}\hline
\mf k & M & {\rm conditions}\\ \hline\hline
\mf{so}(n)& Gr(k,n)&\\ \hline
\mf{spin}(7)& Gr(2,8)& \mf{spin}(7)\subseteq \mf {su}(8)\;
via\;spin\;rep.\\ \hline
\mf{sp}(n)& Gr(k,2n)&  \\  \hline
\mf{sp}(n)& Gr(k,2n+1)& k\neq 3  \\  \hline

%
%
% Caso Punto fisso Irr - Rid
%
%
% Punto fisso caso Irr - Irr
%
\asp (2)+\su (n-4) & Gr (4,n) & n \geq 9 \\ \hline
\su (l) + \su (p) + \su(q)   & Gr (2,n)  & p,q >2,  n=p+q+l
\\ \hline 

\su(l)  + \su (2) + \su (n-l-2)  & Gr (2,n) &  l,n-l>2
\\  \hline
\ \ \ \ \su(l)  +\su (p) + \asp (q)  \ \ \ \ &
\ \ \ \ Gr (2,n) \ \ \ \ \ \  &  \ \ \ q>1,\ p>2, n=p+2q+l \ \\
\hline
%\end{array}
%$$
%$$
%\begin{array}{|c|c|c|}\hline
%\mf k & M & {\rm conditions}\\ \hline\hline
%
\su(l) + \su (n-l-1)  & Gr (k,n) & n >k +l+1, l>k \\ \hline
 \su
(n-1)    & Gr(k,n) &  n \geq 4 \\  \hline
\mf{su}(l)+\mf{su}(n-l) & Gr (k,n)  & k \neq n-l\ {\rm or \ } k
\neq l
\\ \hline
\mf{su}(p)+\mf{su}(q)& Gr (k,n) & q>k, p>k, n=p+q+1 \\ \hline
\end{array}
$$

\end{small}
\end{thm}
All the Lie algebras listed in the first column, unless explicitly specified,
are meant to be standardly embedded into $\su (n),$ e.g.
$\asp (2) + \su (n-4) \subset \su (4)+ \su (n-4) \subset \su (n).$

%
%$\R(\alpha)$ the lines in $\mf z\times \R$ different from
%$y=\alpha x.$ Looking at Table $2,$ one can see that the
%semisimple part of $\mf k$ is at most given by the sum of three
%terms $\mf {k}_1\subseteq \su (l),$ $\mf{k}_2$ and $\mf{k}_3;$ we
%denote by $\R$ the centralizer of $\mf {k}_2+\mf{k}_3$ in $\su
%(n-l),$ finally $\mf{t}_1$ denotes a torus in $\su (2).$ With this
%notation we can also determine exactly the {\em minimal} Lie
%algebras $\mf k$ acting coisotropically on $M.$
%
%
In the following theorem we determine all the
{\em  minimal} non-semisimple
Lie algebras $\mf k$ of the compact connected subgroups of $\SU (n)$
acting
coisotropically on complex Grassmannian.
\begin{thm}\label{nonsemisemplice} Let $K$ be a not semisimple
compact
connected Lie subgroup of
$\SU(n)$ acting non transitively on $M.$ Then $K$ acts
coisotropically on $M$ if and only if its Lie algebra $\mf k,$ up to conjugation,
contains one of the Lie algebras appearing in Table $2,$
and their semisimple parts coincide.
\begin{small}
\begin{center}
{\bf Table 2}
\end{center}
$$
\begin{array}{|c|c|c|}\hline
%
%  1
%
\mf k & M & {\rm conditions}\\ \hline\hline
\mf z +\asp (2)+\su (4) & Gr (4,8) & \\ \hline
\mf z + \su (2) + \asp(n-1)  & Gr (2,2n)   & n \geq 3 \\ \hline
\mf z+ \su (3) + \asp (n-1) & Gr (3,2n+1) & n \geq 3 \\ \hline
\mf z+\mf{su}(k)+\mf{su}(k)& Gr (k,2k) &
\\ \hline\hline
%
%
% 2
%
\mf{t}_1 + \su (n-2) & Gr(2,n) & n>4,\ \mf{t}_1 \mathrm{maximal\ torus \ in}
\ \su(2) \\ \hline
\mf a +\su (n-2)& Gr(k,n)& k > 2,\ n>4 \\ \hline
\mf a + \asp (n) & Gr(3,2n+1)&\\ \hline
\mf{a}+\asp (n-1)    & Gr (2,2n) & n \geq 2 \\ \hline
%    k<= n-l,l
%
\mf a +\su(2 ) & Gr(2,4) & \\ \hline
%
%
% 3
%
%
%\end{array}
%$$
%$$
%\begin{array}{|c|c|c|}\hline
%\mf k & M & {\rm conditions}\\ \hline
\mf a+\su (2) + \su (2) + \su(2)   & Gr (2,6)  &
 \\ \hline
\mf a + \su(2) +  \su (2) + \asp (n-2)  & Gr (2,2n)&  n> 3,\
n=2q+4  \\ \hline
\mf a+ \su(2) + \asp (p) + \asp (q)  & Gr(2,2n) &
p,q \geq 2,\ n=2p+2q+2 \\  \hline
%
%
% 4
%
\mf a + \su (k)  + \su (k-1)   & Gr (k,2k)  &  k\geq 3
\\ \hline
\mf a + \su (k)  + \su (k)   & Gr (k,2k+1)  &  k\geq 2
\\ \hline
\mf a +\su(2) + \mf{sp}(n-1)   & Gr(2,2n+1) & n >2
\\ \hline\hline % qui sotto
\su (n-4) + \R + \su (2) + \su(2)   & Gr (2,n)  &  n> 6
 \\ \hline
\su(l) + \R +  \su (2) + \asp (q)  & Gr (2,n)& l > 2,\ q \geq 2,\
n=2q+2+l  \\ \hline
\su(l) +\R + \asp (p) + \asp (q)  & Gr (2,n) &
p,q \geq 2,\ n=2p+2q+l, l > 2 \\  \hline\hline
\su(l) + \R + \mf{sp}(n-1)   & Gr(2,2n+l-1) & n > 2,\ l > 2
\\ \hline
\su (l) + \R + \su (n-l-1) & Gr (k,n) & n \geq l+3,\ l \geq 2,\ k
>n-l \\ \hline\hline
%
% 5
%
%
\R({\frac{n}{k(n-k)}}) + \su(k)  + \su(n-k-1)   & Gr(k,n) & n >2k+1 \\
\hline
\R ( \frac{ (k-1) n }{n-1} ) + \su (k-1) + \su (n-k) & Gr(k,n) & k
\geq 3 \\ \hline
\R (\frac{qn}{n-1}) + \su (p) + \su (q)& Gr (k,n) & n=p+q+1,\
p>k,\ k \geq q+2. \\ \hline
\R({\frac{-nq}{2(n-2)}}) + \su (2) + \su (p) + \asp (q) & Gr (2,n)
& p>2,\ q>1, \ n=2+p+2q \\  \hline
\end{array}
$$
\end{small}
\end{thm}
$\ $ \\
The notations used in Table $2$ are as follows. We denote with  $\mf{z}$
the one dimensional center of
the Lie algebra $\mf h$ of the group $H=\SUr{k}{n-k},$ with
$\mf a$  a Cartan subalgebra of the centralizer
of the semisimple part of $\mf k.$ In the third block we see that the
semisimple part of $\mf k$ is  given by the sum of three
ideals $\mf {k}_1\subseteq \su (l),$ $\mf{k}_2+\mf{k}_3\subset \su (n-l)$ and
 we
denote by $\R$ the centralizer of $\mf {k}_2+\mf{k}_3$ in $\su
(n-l).$  In the fourth block $\R$ denotes the
centralizer of $\asp (n-1)$ ($\su(n-l-1)$) in $\su (2n-1)$
(resp. $\su (n-l).$
Finally, looking at the last block, we consider in the first three cases
the product $\mf z \times \R$ where $\R$ is the centralizer of the semisimple
part in $\su (n-1)$. With this notation
$\R(\alpha)$ denotes  any  line in $\mf z\times \R$ different from
$y=\alpha x.$ In the last case we use the same notations, pointing out that
here $\R$ is the centralizer of $\su (p)+\asp (q)$ in $\su (p+2q).$\\
\\
Victor Kac \cite{Kac} obtained a complete classification (Tables
Ia, Ib, in Appendix) of irreducible multiplicity free actions
$(\rho,V)$. Most of these include a copy of the scalars $\C$
acting on $V$. We will say that a multiplicity free action
$(\rho,V)$ of a complex group $G$ is {\em decomposable} if we can
write $V$ as the direct sum $V=V_1\oplus V_2$ of proper
$\rho(G)$-invariant subspaces in such a way that
$\rho(G)=\rho_1(G)\times\rho_2(G)$, where $\rho_i$ denotes the
restriction of $\rho$ to $V_i$. If $V$ does not admit such a
decomposition then we say that $(\rho,V)$ is an {\em
indecomposable} multiplicity free action. C. Benson and G.
Ratcliff have given the complete classification of indecomposable
multiplicity free actions (Tables IIa, IIb Appendix). We recall
here their theorem (Theorem 2, pag 154 \cite{BR})
\begin{thm}\label{Ben}
Let $(\rho,V)$ be a regular representation of a connected
semisimple  complex algebraic group $G$ and decompose $V$ as a
direct sum of $\rho(G)$-irreducible subspaces,
$V=V_1\oplus V_2 \oplus \cdots \oplus V_r$.
The action of $(\C^*)^r \times G$
on $V$ is an indecomposable multiplicity free action if and only
if either
\begin{description}
\item[(1)] $r=1$ and $\rho(G)\subseteq \GL(V)$ appears in Table Ia (see the 
Appendix);
\item[(2)] $r=2$ and
$\rho(G)\subseteq \GL (V_1)\times \GL( V_2)$ appears
in Tables IIa and IIb (see the Appendix).
\end{description}
\end{thm}
In \cite{BR} are also given conditions under which one can
{\em remove} or {\em reduce} the copies of the scalars preserving the
multiplicity free action. Obviously  if an action is coisotropic
it continues to be coisotropic also when this action includes
another copy of the scalars. We will call {\em minimal} those
coisotropic actions in which the scalars, if they appear, cannot
be reduced.\\
Let $K$ be a compact group acting isometrically on a
compact K\"ahler manifold $M.$ We say that $M$ is
$K^\C$-{\em almost homogeneous} if $K^\C$ has an open orbit in $M.$ If all
Borel subgroups of $K^\C$ act with an open orbit on $M$, then the
$K^\C$-open orbit $\Omega$ is called a
{\em spherical homogeneous space} and $M$ is called a
{\em spherical embedding} of $\Omega.$
We will briefly recall some results that will be used in the
sequel.
\begin{thm}{\label{class}}\cite{Hw}
Let $M$ be a connected compact K\"ahler manifold with an isometric
action of a connected compact group $K$ that is also Poisson. Then
the following conditions are equivalent:
\begin{itemize}
\item [(i)] The $K$-action is coisotropic.
\item [(ii)] The cohomogeneity of the $K$ action is equal to the
difference between the rank of $K$ and the rank of a regular
isotropy subgroup of $K$.
\item [(iii)] The moment map $\mu\,:\,M\rightarrow{\mathfrak{k}^*}$
separates orbits.
\item [(iv)] The K\"ahler manifold $M$ is projective algebraic,
$K^\mathbb{C}$-almost homogeneous and a spherical embedding of the
open $K^\mathbb{C}$-orbit.
\end{itemize}
\end{thm}
\n
As an immediate consequence of the above theorem one can deduce,
under the same hypotheses on $K$ and $M,$ two simple facts that
will be frequently used in our classification:
\begin{itemize}
\item [1] Let $p$ be a fixed point on $M$ for the $K$-action, or $Kp$
a complex $K$-orbit, then the $K$-action is coisotropic if and
only if the slice representation is coisotropic.
\item [2] {\em dimensional condition.} If $K$ acts coisotropically on $M$
the dimension of a Borel subgroup $B$ of $K^\C$ is not less than
the dimension of $M$.
\end{itemize}
A relatively large class of coisotropic actions is provided by
polar ones. We recall here that an isometric action of a group $K$
is called {\em polar} on $M$ if there exists a properly embedded
submanifold $\Sigma$ which meets every $K$-orbit and is orthogonal
to the $K$-orbits in all common points. Such a submanifold is
called a {\em section}, and, if it is flat, the action is called
{\em hyperpolar}.\\
A result due to Hermann states that given $K$  a compact Lie group
and two symmetric subgroups $H_1$,$H_2\subseteq{K}$,  then $H_i$
acts hyperpolarly on $K / H_j $ for $i,{j}\in{1,2}.$ Therefore, given
an Hermitian symmetric space as $Gr(k,n)= \SU(n)/ \SUr{k}{n-k}$ it
is clear that subgroups like $\SO(n)$, $\Sp(n)$ and  $
\SUr{l}{n-l}$ will act on $Gr (k,n)$ hyperpolarly, hence, since for
\cite{PoT}
a  polar action on an irreducible compact
homogeneous K\"ahler manifold is co\-iso\-tro\-pic, these actions are
coisotropic.\\
Once we have determined the complete list of coisotropic actions
on $Gr(k,n),$ in Table $2,$  we have also investigated which ones
are polar. Dadok \cite{Da}, Heintze and Eschemburg \cite{He} have
classified the irreducible polar linear representations, while
I.Bergmann \cite{Bergmann} has found all the reducible polar ones.
Using their results we determine in section $3$
  the  polar actions on complex Grassmannians. We get the following
\begin{thm}
 Let $K$ be a compact
connected Lie subgroup of $\SU(n)$
acting non transitively on the complex Grassmannians
$M=Gr(k,n)= \SU(n)/ \SUr{k}{n-k}.$  Then $K$ acts polarly
on $M$ if and only if its Lie algebra $\mf{k}$ is  conjugate
to one of the Lie algebras  appearing in Table $3.$
\begin{small}
\begin{center}
{\bf Table 3}
\end{center}
$$
\begin{array}{|c|c|c|}\hline
\mf k & M & {\rm conditions}\\ \hline\hline
\mf{so}(n)& Gr(k,n)&\\ \hline
\mf{sp}(n)& Gr(k,2n)&  \\  \hline
\R+ \mf{su}(l)+\mf{su}(n-l)& Gr (k,n) &  \\ \hline
\mf{su}(l)+\mf{su}(n-l)& Gr (k,n) & l\neq n-l \\ \hline
\end{array}
$$
\end{small}
\\
where $\R$ is the centralizer of $\su (n-l)+ \su (l)$ in $\su (n)$.
In particular the $K$-action is hyperpolar.\label{polar}
\end{thm}
We here briefly explain our method in order to prove our main theorem.
Thanks to Theorem \ref{class},(iv) we have that
 if $K$ is a subgroup of a compact Lie
group $L$ such that $K$ acts coisotropically on $M$ so does $L$.  As a consequence, in order to
classify coisotropic actions on  $Gr(k,n),$ one  may suggest a
sort of ``telescopic'' procedure by restricting  to maximal
subgroups $K$ of $\SU(n)$,
 hence passing to
maximal subgroups that give rise to coisotropic actions and so
on.\\ In the Appendix  we give the complete
list of maximal subgroups $K$  of $\SU(n)$.\\
This paper is organized as follows. In section 2 we prove Theorem
\ref{risultato} and  Theorem \ref{nonsemisemplice}. We have
divided this section in five subsections in each of which we
analyze separately one of the maximal subgroups of
 $\SU(n)$.
In the third section we give the proof of Theorem \ref{polar}.\\
We enclose, in the Appendix, the tables of irreducible  and
reducible linear multiplicity free representations (Tables Ia, Ib
and Tables IIa, IIb respectively), the table of maximal subgroups
of $\SU(n)$ and $\Sp(n)$ (Tables III, IV) and the table of reduced
prehomogeneous triplets (Table V), since  we will frequently refer to them in our analysis.
%%%%%%%%%%%%%%%%%%%%%%%%%%%%%%%%%%%%%%%%%%%%%
%%%%%%%%%%%%%%%%%%%%%%%%%%%%%%%%%%%%%%%%%%%%%
%%%%%%%%%%%%%%%%%%%%%%%%%%%%%%%%
%
%
\section{Proof of Theorem \ref{risultato} and Theorem \ref{nonsemisemplice}}
In the following subsections we will go through all maximal
subgroups $K$ of $\SU(n)$ according to Table III in Appendix.
%
%
%%%%%%%%%%%%%%%%%%%%%%%%%%%%%%%%%%%%%%%%%%%%%%%%%%%%%%%%%%%%%%%%%%%%%%%%%%%%
%                         L'AZIONE DI SO(n)
%
%%%%%%%%%%%%%%%%%%%%%%%%%%%%%%%%%%%%%%%%%%%%%%%%%%%%%%%%%%%%%%%%%%%%%%%%%%%%%%%%%%%%%%%%%%%%%%%%%%%%%%%%%%%%%%%%%%%%%%%%%%%%
\subsection{The case $K=\SO(n)$}
Note that $\SO(n)$ is a symmetric subgroup of $\SU(n)$. The action
is polar hence coisotropic on complex Grassmannians. We can now
investigate the behaviour of subgroups of $\SO(n)$. Take a
$k$-plane $\pi$ in $\mathbb{R}^n.$ The $\SO(n)$-orbit of $\pi$ is
a real non oriented Grassmannian. Note that the orbit coincides
with the set of fixed points of the conjugation, which is an
antiholomorphic isometry on $Gr (k,n)$, hence the orbit is totally
real in $Gr (k,n)$. The following proposition is easy to prove
using well known properties of the moment map.
\begin{prop}Let $K$ be a semisimple Lie group which acts on a
K\"ahler manifold $M$ with a totally real orbit $Kp$ and moment
map $\mu$. Then $\mu(Kp)=0.$
\end{prop}
In the sequel we will identify $\mathfrak{k}\cong{\mathfrak{k}^*}$
by means of an $Ad(K)$-invariant inner product  on $\mathfrak{k}.$
\\
From the previous proposition we deduce that $0$ belongs to the
image of the moment map
$\mu\,:Gr(k,n)\rightarrow\,\mathfrak{so(n)}.$ Let now take the
subgroup $K'$ of $\SO(n)$. Obviously its image via the moment map
$\mu'$ of $K',$ given by the composition of $\mu$ with the
projection map $\mf{k}\rightarrow \mf k'$, contains $0.$ If $K'$
acts coisotropically, then $\mu'$ separates, by Theorem
\ref{class} (iii), the orbits and therefore ${\mu'}^{-1}(0)$
coincides with the real non oriented Grassmannian. The problem is
then restricted to the analysis of those subgroups $K'$ of
$\SO(n)$ that act transitively on the real Grassmannian and
coisotropically on the complex $Gr (k,n).$ The first problem was
solved by Onishchick in \cite{On}.\\ We have to distinguish two
cases, namely $k=2$ and $k>{2}.$\\ In the first case we have to
consider the subgroups $G_2$ and $Spin(7)$ of $\SO(7)$ and
$\SO(8)$ respectively. The first one, $G_2,$ cannot act
coisotropically on $Gr (2,7)$ because the Borel subgroup of
$G_2^\C$ has complex dimension $8< \dim {Gr(2,7)}$. $Spin(7)$
acts, instead, coisotropically on $Gr (2,8).$ Note that the
totally real orbit is, in this case, the non oriented real
Grassmannian $Gr_\R (2,8)$ and $\U(3)$ is the isotropy subgroup of
$Spin(7)$. It is well known that the cohomogeneity
$\textrm{chm}(Gr(2,8),Spin(7))=\textrm{chm(Slice},\U(3))$; now,
the slice representation of $\U(3)$ on $Gr(2,8)$ can be deduced
immediately from the isotropy representation of $\U(3)$, that is
given by $\C^3\oplus (\C^3)^*=\C^3\oplus \Lambda^2(\C^3)$, since
the orbit is totally real.  Furthermore observe that
$\textrm{rk}(Spin(7))=\textrm{rk}(\U(3))$ and the same holds for
the principal isotropy subgroups; hence, by Theorem \ref{class}
(ii), the $Spin(7)$-action is coisotropic if and only if the
$\U(3)$-action is coisotropic. Now, looking at the list of Benson
and Ratcliff (Table IIa, Appendix), we conclude that the action of
$\U(3)$ is coisotropic since the scalars act on $\C^3\oplus
(\C^3)^*$ as $z(v_1,v_2\wedge v_3)=[zv_1,z^2(v_2\wedge v_3)]$ and
$a=1\neq -b=-2$. \\ In the second case, when $k>2$, only
$K'=Spin(7)$ acts transitively on the real non oriented
Grassmannian of 3-planes in $\mathbb{R}^8.$ As for $G_2$, the
Borel subgroup of $Spin(7)$ is too small: it has, in fact, complex
dimension 12 which is smaller than the complex dimension of $Gr
(3,8).$\\ This completes the analysis ; we have thus proved that
the only subgroups of $\SO(n)$ which act coisotropically on some Grassmannian
 $Gr(k,n)$ are $\SO(n)$ and  $Spin(7)\subset \SO(8)$. 
%%%%%%%%%%%%%%%%%%%%%%%%%%%%%%%%%%%%%%%%%%%%%%%%%%%%%%%%%%%%%
%%%%%%%%%%%%%%%%%%%%%%%%%%%%%%%%%%%%%%%%%%%%%%%%%%%%%%%%%%%%%
%                     L'AZIONE DI SU(p)\times SU(q)

%%%%%%%%%%%%%%%%%%%%%%%%%%%%%%%%%%%%%%%%%%%%%%%%%%%%%%%%%%%%%%%%%%%%%%%%%%%%%%%%%%%%%%%%%%%%%%%%%%%%%%%%%%%%%%%%%%%%%%%%%%%%
\subsection{The  case $K=\SU(p) \otimes \SU(q)$.}
For simplicity we will denote the group  $\SU(p) \times \SU(q)$ by
$K\subseteq \SU(n)$, with $n=pq.$\\ In this case we will determine
the slice representation on a complex orbit, showing that is not
coisotropic. We recall the following simple lemma
\begin{lem}\label{pino}Let $\rho:G \rightarrow \GL (V)$ be a linear
coisotropic
representation. If $V$ splits as the direct sum $V_1 \oplus V_2
\oplus \ldots \oplus V_l$
of irreducible submodules then $\rho(G)$ restricted on each $V_i$ is
coisotropic.\end{lem}
We have to distinguish two different cases.
\begin{enumerate}
\item {\em Either  $p$ or $q$ is bigger than $k$}.
We can assume, for example, that $p \geq k$. We choose a
$k$-dimensional subspace $W\subseteq \C^p$ and a one-dimensional
$V\subseteq \C^q.$  The stabilizer $K_\pi$, where $\pi$ is the
$k$-plane $W \otimes V \subseteq \C^p \otimes \C^q,$ is $
\SUr{k}{p-k}\times \SUr{1}{q-1}$; it is easy to see that the tangent
space $T_\pi{G_r(k,p)}=\pi^*\otimes \pi^\perp$ splits, as a
$K_\pi$-module, as
\begin{small}
$$
(W^* \otimes W^\perp) \oplus (V^* \otimes V^\perp) \oplus
(\mf{sl}(W)\otimes V^* \otimes V^\perp) \oplus (W^*\otimes
W^\perp\otimes V^*\otimes V)$$
\end{small}
where $W^\perp,V^\perp$ are the orthogonal complements of $W,V$ in
$\C^p,\C^q$ respectively. Here we have used the fact that
$W^*\otimes W=\C\oplus \mf{sl}(W)$ as $U(k)$-modules. The first
two factors can be identified with the tangent space at $\pi$ of
the $K$-orbit $K\pi,$ which is therefore complex; on the other
hand the summand $\mf{sl}(W)\otimes V^* \otimes V^\perp$ does not
appear in Kac's list, so that the $K$-action is not coisotropic, thanks to Lemma \ref{pino}.
\n
\item {\em $p$ and $q$ are smaller than $k$.}
Clearly $ k=mp \ + \ l$ where $m<q$ and $l<p$. Choose an
$l$-dimensional subspace $W\subseteq \C^p,$ an $m$-dimensional
$V\subseteq \C^q,$ and $ v_{m+1}\in \C^q \setminus V.$ We consider
a $k$-complex plane $\pi=(\C^p \otimes V) \oplus (W \otimes
v_{m+1})$; it is easy to check that the $K$-orbit $K \pi\cong Gr(l,p)
\times \SU (q) / \Su (\U (1) \times \U (m) \times \U (q-m-1))$ is
complex. Following the same procedure used in the previous case,
with slightly heavier computations, we determine the slice summand
$\mf{sl}(W^{\perp}) \otimes V^* \otimes v_{m+1}$, on which
$Ad(\U(l)) \otimes \SUr{m}{1} $ acts; this action does not appear
in Kac's list. Hence, using again Lemma \ref{pino}, we can conclude that
the $\SU(p) \otimes \SU(q)$-action cannot be coisotropic.
\end{enumerate}
%%%%%%%%%%%%%%%%%%%%%%%%%%%%%%%%%%%%%%%%%%%%%%%%%%%%%%%%%%%%%%%%%%%%%%%%%%%%%%%%%%%%%%%%%%%%%%%%%%%%%%%%%%%%%%%%%%%%%%%%%%%%
%
%                       L'AZIONE DI UN GRUPPO SEMPLICE
%
%%%%%%%%%%%%%%%%%%%%%%%%%%%%%%%%%%%%%%%%%%%%%%%%%%%%%%%%%%%%%%%%%%%%%%%%%%%%%%%%%%%%%%%%%%%%%%%%%%%%%%%%%%%%%%%%%%%%%%%%%%%%
\subsection{The case  of a simple Lie group $K=H$ such that  $\rho(H)$
is an irreducible representation  of  complex type.}
This case can be analyzed using the work of Sato and Kimura \cite{SK}.
 Throughout this section we will identify the fundamental highest weights $\Lambda_l$ with the corresponding irreducible representations.\\ We recall
here that, given $G$ a connected complex linear algebraic group,
and $\rho$ a rational representation of $G$ on a finite
dimensional complex vector space $V,$ such a triplet $(G,\rho,V)$
is {\em prehomogeneous} if $V$ has a Zariski dense $G$-orbit.
\begin{defi} Two triplets $(G,\rho,V)$ and $(G',\rho',V')$ are called
{\em equivalent} if there exist a rational isomorphism
$\sigma:\rho(G)\to\rho'(G')$ and an isomorphism $\tau:V\to V'$,
both defined over $\C$ such that the diagram is commutative for
all $g\in G$.
$$\xymatrix{
V\ar[r]^-\tau \ar[d]_-{\rho(g)} & V' \ar[d]^-{\sigma(\rho(g))} \\
V\ar[r]^-\tau    &   V'}$$ This equivalence relation will be
denoted by $(G,\rho,V)\cong (G',\rho',V')$.
\end{defi}
 We say
that two triplets $(G,\rho,V)$ and $(G',\rho',V')$ are {\em
castling transforms} of each other when there exist a triplet
$(\tilde{G},\tilde{\rho},V(m))$ and a positive  number $n$ with
$m>n\geq 1$ such that $$(G,\rho,V)\cong (\tilde{G}\times \SL(n)
,\tilde{\rho}\otimes \Lambda_1,V(m)\otimes V(n))$$ and
$$(G',\rho',V')\cong (\tilde{G}\times \SL(m-n)
,{\tilde{\rho}}^*\otimes \Lambda_1,V(m)^*\otimes V(m-n)),$$ where
${\tilde{\rho}}^*$ is the contragradient representation of
${\tilde{\rho}}$ on the dual vector space $V(m)^*$ of $V(m)$. We
recall that $V(n)$ is a complex vector space of dimension $n.$ A
triplet $(G,\rho,V)$ is called {\em reduced} if there is no
castling transform $(G',\rho',V')$ with $\dim V'<\dim V.$ We give
here a useful lemma that permits to find relations between
coisotropic actions on Grassmannians,  which must have an open
dense orbit, and prehomogeneous triplets. This can be easily
proved as a consequence of
Lemma 5 and Proposition 7 pag. 37 in section 3 \cite{SK}.
\begin{lem}
Let $G$ be any complex, connected Lie group. $G$ acts with an open
dense orbit on $Gr (k,n)$ if and only if $G \times GL(k)$ acts
with an open dense orbit on $\C^n \otimes \C^k$ i.e. $(G\times \GL
(k),\rho,\mathbb{C}^n\otimes \mathbb{C}^k)$ is a prehomogeneous
triplet.
\end{lem}
Let $H$ be a simple group; denote  with $d$ the degree of $\rho$
and with $V(d)$ its  representation space. We are looking for
simple groups $H$, whose representation is irreducible and of
complex type, that act coisotropically on complex Grassmannians
$Gr(k,d)$. Recall that one can restrict the analysis to
$2<k\leq{\frac{d}{2}}$; in fact for $k>\frac{d}{2}$ one can take
the dual Grassmannian $Gr(d-k,d).$ With these notations we state
and prove the following
\begin{lem} If $\mf h \neq A_n$ then the triplet
$(H\times \GL(k), \rho\otimes \Lambda_1,V(d)\otimes V(k))$ is reduced.
\end{lem}
\begin{proof}
Suppose, by contradiction, that the triplet is not reduced. Then
there would exist two triplets $(\tilde{G},\tilde{\rho},V(m))$ and
$(G',\rho',V'),$ with $\dim V' < \dim (V(d)\otimes V(k)),$ such
that
\begin{description}
\item[(a)] $(H \times \GL (k),\rho \otimes
\Lambda_1,V(d) \otimes V(k)) \cong (\tilde{G}\times \SL(n),
\tilde{\rho}\otimes \Lambda_1,V(m)\otimes V(n));$
\item[(b)]  $(G',\rho',V')\cong (\tilde{G}\times \SL(m-n)
,{\tilde{\rho}}^*\otimes \Lambda_1,V(m)^*\otimes V(m-n)).$
\end{description}
Switching to the Lie algebra level we get, from the first
condition $\mf h + \mf{gl} (k) \cong d \tilde{\rho}
(\tilde{\mf{g}}) + \mf{sl} (n)$ which implies, since $\mf h \neq
A_n,$ that $n=k$ and $m=d$, and, from the second, we get $\dim
V'=d(d-k)$. This gives rise to a contradiction since $\dim V' \geq
\dim (V(d)\otimes V(k))$
\end{proof}
\n Sato and Kimura have completely classified the prehomogeneous
reduced triplets. Hence the problem, if $\mf h\neq A_n$, boils down to see whether the triplet $(H \times \GL (k),\rho \otimes
\Lambda_1,V(d) \otimes V(k))$ appears in the list of Sato Kimura
(see Table V Appendix). If it is not the case the $H$-action is
not
coisotropic.\\
We find, for $k=2,$ the triplet $(Spin(7,\C)\times \GL (2),
\rm{spin. \ representation}\otimes \Lambda_1,\C^8\otimes \C^2)$;
hence $Spin(7)$ has an open orbit on $Gr (2,8);$ it has already
been studied as subgroup of $\SO(8).$ This action is coisotropic.
\\ The other groups, arising from the triplets $(Spin(10)\times
\GL (2), \rm{half-spin. \ representation}\otimes \Lambda_1,
\C^{16}\otimes\C^2)$, $(G_2\times \GL (2), \Lambda_2\otimes
\Lambda_1, \C^7\otimes \C^2)$ and $(E_8^\C \times \GL
(2),\Lambda_1\otimes \Lambda_1, \C^{27}\otimes \C^2)$ can be
excluded thanks to the {\em dimensional condition}.\\ For $k>2$ we
find again $Spin(7)$ on $Gr (3,8)$, and $Spin(10)$ on $Gr (3,16)$
and they are excluded for the same dimensional reason.\\ Note that
in the list of Sato and Kimura the group $\Sp(n)$ appears; this
case will be studied
separately in next paragraph.\\
\\
Now, suppose $\mathfrak{h}=A_n$. We impose the {\em dimensional
condition}. Denote with $\mathfrak{b}$ the Lie algebra of the
Borel subgroup of $H^\C$. We get the following inequality
$$\dim_{\mathbb{C}} \mathfrak{b} \geq \dim_{\mathbb{C}}Gr(k,d)$$
that, when $\mathfrak{h}=A_n,$ becomes
\begin{equation}
\label{1}\frac{(n-1)(n+2)}{2k}+k\geq d.
\end{equation}
Using Young diagrams we can write down the degree of each
representation. \\To solve this case we distinguish the case $k=2
$ and $k\neq 2.$ For $n=2$ and $k=2$ the fact that only
$\Lambda_1$ and $2 \Lambda_1$ are admissible is straightforward.
If $k=2,$ and $n \geq 3,$ then it is easy to prove that $d<
\frac{1}{2}n(n+1) < n^2$ hence, using
 Proposition 7 and  Corollary 8 (pag. 45 section 3 of \cite{SK}), we
have only to investigate the  representations $\Lambda_1,\
\Lambda_2,\ \Lambda_{n-1}, \Lambda_{n-2}, \ n=3,4.$  Note that
$\Lambda_1, \Lambda_{n-1}$ correspond to transitive action of $\SU
(n)$  on $Gr(k,n).$ We then find the following two triplets: when
$n=2,$ $(\SL (2), 2 \Lambda_1, S^2 ( \C^2))$ which corresponds to
a coisotropic action on the complex projective space and for
$n=4,$ we get $(\SL (4),  \Lambda_2, \C^6 )$ that must be excluded
since it
is of real type.\\
In the general case, when $2<k\leq \frac{d}{2},$  we can do
analogous calculations excluding all the simple groups which act
non transitively, i.e which have representations different from
$\Lambda_1$ and $\Lambda_{n-1}.$ Indeed we note that
$$
d\leq \frac{(n-1)(n+2)}{2k} + k\leq \frac{(n-1)(n+2)}{6}+
\frac{d}{2}
$$
hence
\begin{equation} \label{2}
d \leq \frac{(n-1)(n+2)}{3}.
\end{equation}
Consider the case $\rho = \Lambda_2 $ separately. The inequality
(\ref{2}) holds for $n\leq 4$ and the corresponding cases are
 $(\SL (3),\Lambda_2,\mathbb{C}^3),$ which corresponds to the
transitive action of $\SL (3)$ on $Gr (1,3),$ and to $(\SL
(4),\Lambda_2,\mathbb{C}^6)$ that can be excluded since it is of
real type.\\ 
We then study the general case, analyzing $d(\Lambda_l)$, where $2<l\leq \frac{n}{2}$
(if $l>\frac{n}{2}$ we can consider the dual representation which
shares
 the same behavior). We shall prove that the inequality (\ref{2})  does
not hold if $\rho \neq \Lambda_1, \Lambda_{n-1}.$ Indeed it
becomes for $\rho=\Lambda_l$, $$3n(n-2) \cdots
(n-l+1)\leq(n+2)l!$$ with the condition $n\geq 2l>4.$ This can
never be satisfied.\\ Now, observing that the degree $d(\sum _i
m_i\Lambda_i)$ increases as $m_i$ increases it is sufficient to
prove that $m\Lambda_1, \Lambda_1+\Lambda_2$ and $\Lambda_1+
\Lambda_{n-1}$ do not satisfy condition (\ref{2}) in order to get
our assert. This is a straightforward calculation using Young
diagrams.\\ We conclude that no simple group $H$ with $\rho(H)\neq
\Lambda_1,\Lambda_{n-1}$ gives rise to a coisotropic action, except for
$Spin(7)\subset \SO(8)$ and $\Sp(n)\subset \SU(2n)$.
%%%%%%%%%%%%%%%%%%%%%%%%%%%%%%%%%%%%%%%%%%%%%%%%%%%%%%%%%%%%%%%%%%%%%%%%%%%%%%%%%%%%%%%%%%%%%%%%%%%%%%%%%%%%%%%%%%%%%%%%%%%%
%
%                            L'AZIONE di Sp(n)
%
%%%%%%%%%%%%%%%%%%%%%%%%%%%%%%%%%%%%%%%%%%%%%%%%%%%%%%%%%%%%%%%%%%%%%%%%%%%%%%%%%%%%%%%%%%%%%%%%%%%%%%%%%%%%%%%%%%%%%%%%%%%%
\subsection{ The case $K=\Sp(n)$.} We now consider the case
$K=\Sp(n)\subseteq \SU(2n)$ acting on the Grassmannian $Gr
(k,2n)$. Note that $\Sp(n)$ is a symmetric subgroup of $\SU(2n)$;
the action is then  polar hence coisotropic on complex
Grassmannians. We will now go through the maximal subgroups $K'$
of $Sp(n)$ (see Table IV in the Appendix, for the complete list).\\
\\
\noindent (i) $K'=\U(n).$ We  recall that the subgroup $K'=\U(n)$
acts on $\C^{2n}$ reducibly, namely as $\rho_n \oplus \rho_n^*$
where $\rho_n$ denotes the standard representation of $U(n)$ in
$\C^n$. We consider a $k$-dimensional subspace $\pi \subseteq
\C^n$, since we can always suppose $k\leq n$; then $K' \pi = \U (n)
/ \U (k) \times \U (n-k)$ and $T_\pi Gr (k,2n)$ splits, as a
$\U (k) \times \U (n-k)$-module, as $ S^2 (\pi^*) \oplus \Lambda^2
(\pi^*) \oplus \pi^* \otimes \C^{n-k} \oplus \pi^* \otimes
{\C^{n-k}}^*.$ This means that $K'\pi$ is a complex orbit with
slice representation $S^2 (\pi^*) \oplus \Lambda^2 (\pi^*) \oplus
\pi^* \otimes {\C^{n-k}}^*$ which is not coisotropic by \cite{BR}.
\\
\\
\n
(ii) $K'=\SO(p) \otimes \Sp(q),\ pq=n, p \geq 3,\ q \geq 1.$
Clearly $$SO(p) \otimes \Sp(q) \subset \SU(p) \otimes \SU(2q)$$
hence the action cannot  be coisotropic as the action of $\SU(p)
\otimes \SU(2q)$ is not (section 2.2).\\
\\
\n (iii) $K'=\rho(H), H \ simple, \rho \in Irr_{\H}, deg\rho =n.$
With the same estimates on the degree of the representation that
we have used in the previous paragraph, it is easy to prove that
these groups cannot acts coisotropically on the complex
Grassmannians.\\
\\
\n
(iv) $K'=\Sp(l) \times \Sp(n-l).$ First of all we shall  give a
complex orbit for the  $\Sp(n)$-action. We recall that $k \leq
n.$ We identify the $\Sp(n)$-module $\H^n$ with $\C^{2n}=\C^n
\oplus \C^n$ endowed with the anti-holomorphic map $J\in \C^{2n}$
commuting with the $\Sp(n)$-action. If we select a $k$-plane
$\pi\subseteq \C^n$ it is not difficult to see that the stabilizer
$\Sp(n)_\pi$ is given by $\U (k) \times \Sp(n-k)$. We have that
$\pi^{\perp}=J\pi \oplus W,$ where $W$ is a quaternionic subspace
of complex dimension $2(n-k)$; note that $J\pi=\pi^*$ as
$\U(k)$-modules. The tangent space $T_\pi{Gr(k,2n)}$ splits into
irreducible complex $\Sp(n)_\pi$-submodules as
$$\pi^* \otimes \pi^{\perp}=
S^2{\pi^*}\oplus \Lambda^2{\pi^*}\oplus (\pi^*\otimes W),
$$
while on the other hand the isotropy representation of
$\Sp(n)/ \Sp(n)_\pi$ is given by $S^2 \pi^*\oplus L$, where
$L=\H^k\otimes_{\H } \H^{n-k}$ as real submodules. From this we see
that the slice is isomorphic to $\Lambda^2 \pi^*$ and that the
orbit $\Sp(n)\pi$ is complex. Note that this action appears in
\cite{Kac}, this is another way to see that $\Sp(n)$ acts on
complex Grassmannians coisotropically.\\
We now consider $K'=\Sp(l)\times \Sp(n-l)$. We take $k_1, \ k_2$
such that $k_1 \leq l, \ k_2 \leq n-l$ and $k=k_1 + k_2$. Let
$\pi_1$ and $\pi_2$ be two complex $k_1$- and $k_2$-planes in
$\C^{2l}$ and $\C^{2n-2l}$ respectively, such that the orbits in
$Gr (k_1,2l)$ and $Gr (k_2,2n-2l)$ resp. are complex. Let $\pi=
\pi_1 \oplus \pi_2$ be the $k$-complex plane in $\C^{2n}$. It is
easy to see that the stabilizer $K'_\pi$ of the $k$-plane $\pi$ is
 $\U(k_1) \times \Sp( l- k_1) \times \U(k_2) \times \Sp( n-l - k_2)$.
With the same notation that we have used before, it is easy to
check that the slice representation is given by
$$ \Lambda^2
(\pi_1) \oplus \Lambda^2 (\pi_2) \oplus  (\pi_1^* \otimes \pi_2
^* )\oplus  (\pi_1^* \otimes \pi_2^* ) 
\oplus (\pi_1^* \otimes W_2) \oplus  (\pi_2^* \otimes W_1) $$
and looking at Tables IIa, IIb  we can conclude that this action
is not coisotropic on the complex Grassmannians.\\ Note that if we
consider the $\Sp(n)$-action on the complex Grassmannians, with
$k$  even, the orbit of  the $k$-plane $\pi$ generated by
$v_1,v_2\cdots v_k$ and $Jv_1,Jv_2\cdots Jv_k$, is totally real
(for $k=2$ it is $\mathbb{HP}^{n-1}={\Sp(n)}/ \Sp(1) \times \Sp
(n-1)$). We may restrict our analysis, as in the $\SO(n)$ case, to
those subgroups of $\Sp (n)$ that act transitively on the totally
real orbit, and coisotropically on the complex Grassmannians. In
\cite{On} one can find that no subgroups of $\Sp(n)$ can act
transitively on  the totally real orbit. This gives another proof,
when $k$ is even, of the fact that no subgroup of $\Sp(n)$ can act
coisotropically on $Gr (k,n)$.
%%%%%%%%%%%%%%%%%%%%%%%%%%%%%%%%%%%%%%%%%%%%%%%%%%%%%%%%%%%%%%%%%%%%%%%%%%%
%%%%%%%%%%%%%%%%%%%%%%%%%%%%%%%%%%%%%%%%%%%%%
%
%                    IL CASO RIDUCIBILE CON k=l
%
%%%%%%%%%%%%%%%%%%%%%%%%%%%%%%%%%%%%%%%%%%%%%%%%%%%%%%%%%%%%%%%%%%%%%%%%
%%%%%%%%%%%%%%%%%%%%%%%%%%%%%%%%%%%%%%%%%%%%%%%%%%
\subsection{The reducible case}
We can investigate directly the subgroups $L$ of $K=\SUr{l}{n-l}$
because $K$ acts coisotropically on the complex Grassmannians.\\
 Recall that,
in the sequel, we always assume $1<k\leq \frac{n}{2}$. If
$k>\frac{n}{2}$ we refer to the dual Grassmannian $Gr(k',n),$
where $k'=n-k\leq \frac{n}{2}.$\\ With this assumption we divide
the reducible case in three paragraphs: in the first one we will
analyze the case $k=l$, namely when the group $K$ has a fixed
point in $Gr (k,n)$; in the second we suppose $k<l,n-l$, finally
in the third one we complete the analysis studying the case
$k<l,\;k>n-l$. We will also, for simplicity switch to the Lie
algebra level, denoting by $\mf z$ the one dimensional center of
the Lie algebra $\mf k$ of $K.$
%
%
%
%
%
%                        CASO  PUNTO  FISSO
%
%
%
\subsubsection{The fixed point case}
 Consider
$L\subseteq K=\SUr{k}{n-k}$ and $\pi= \C^k$ a fixed $k$-plane. Let
$\mf l$  be the Lie algebra of $L,$ $\mf l \subseteq \mf{su} (k) +
\mf{su} (n-k) + \R= \mf k.$ Suppose $\mathfrak{l}$ acts
coisotropically. We consider the projections $\rho_1 :\mf l
\longrightarrow \mf{su} (k),$ $\rho_2 :\mf l \longrightarrow
\mf{su} (n-k)$ and we put $\mf{l}_i= \rho_i (\mf l ).$ This means
that $\mf z + \mf{l}_1 + \mf{l}_2$ acts coisotropically. We
distinguish some cases according to the irreducibility of $\mf
{l}_i$ on the slice $S$ at $\pi = \C^k,$ that  is given by $(\C^k
)^* \otimes \C^{n-k}.$ In the sequel we refer to Tables Ib, IIa and
IIb in the Appendix, for all the conditions under which one can
remove the scalars preserving the multiplicity free action.\\
\\
{\bf(1)} $\mf{l}_1$ and $\mf{l}_2$ act both irreducibly on $( \C^k
)^*$ and $\C^{n-k}$ respectively. Since $\mf{l}_i \subseteq
\mf{su}(k)$ or $\mf{su}(n-k),$ we have that $\mf{l}_i$ are
centerless. Hence $\mf{l}_1 + \mf{l}_2$ is semisimple and acts
irreducibly  on $(\C^k)^* \otimes \C^{n-k}$ from Table Ia we get
\begin{center}
\begin{tabular}{|c|c|c|c|}
n. & $\mf{l}_1$   & $\mf{l}_2$      & ${\rm conditions}$ \\ \hline
1. & $\mf{su} (k)$ &  $\mf{su} (n-k)$   &   \\
2. & $\mf{sp} (2)$   & $\mf{su} (k)$  &  $ k \geq 4$  \\
3. & $\mf{su} (2)$ & $\mf{sp} (\frac{n-2}{2} )$ & $n \geq 2$ \\
4. & $\mf{su} (3)$ & $\mf{sp} (\frac{n-3}{2} )$ & \\
\end{tabular}
\end{center}
Note that each $\mf{l}_i$ has trivial centralizer in $\mf{su}(k)$
(or in $\mf{su} (n-k)$), so that the centralizer of $\mf{l}$ in
$\mf k$ coincides with the center $\mf{z}.$ Moreover, except in
case 1, when $k=n-k,$ we have that $\mf{l}_1$ is not isomorphic to
$\mf{l}_2$, so that $\h \cong \li 1 +\li 2$ or $\h \cong \li 1 +
\li 2 + \mf z$ or $\h \cong \su (\frac{n}{2})^{\Delta} \subset \su
(\frac{n}{2}) + \su (\frac{n}{2}) +\mf z$ or $\h \cong \su
(\frac{n}{2})^{\Delta} + \mf z.$ Note also that the scalars are
always needed except in  case $1$, $k \neq \frac{n}{2}$, and case
$2$, $k \geq 5$, while the center $\mf z$ always acts
non-trivially on the tangent space. The case $\h \cong \mf z+\su
(\frac{n}{2})^{\Delta}$ can be excluded simply by the dimensional
condition. Summing up we have the following minimal coisotropic
actions
$$
\begin{array}{|c|c|c|}
\h &   M     &  {\rm conditions} \\ \hline
\mf{su} (k) + \su (n-k) & Gr(k,n) & 2k \neq n \\
\mf z + \mf{su}(k) + \su (k)  & Gr (k,2k) &\\
\asp (2)+\su (n-4) & Gr (4,n) & n \geq 9 \\
\mf z + \asp (2) +\su(4) & Gr(4,8) & \\
\mf z  + \su (2) + \asp(n-1)  & Gr(2,2n) & n \geq 3 \\
\mf z  +  \su (3) + \asp (n-1) & Gr (3,2n+1)& n \geq 3\\
\end{array}
$$
\\
{\bf (2)} Suppose $\li{i}$ acts irreducibly and $\li j$ acts
reducibly ($i \neq j, 1 \leq i,j \leq 2$). Assume, for example,
$\li 1$ acts irreducibly on $(\C^k )^*$ and $\li 2$ reducibly on
${\C^{n-k}}.$ Thus the second factor splits as the direct sum $V_o
\oplus  V^1 \oplus \cdots \oplus V^j \oplus V_1 \oplus \cdots
\oplus V_r$ of irreducible $\li 2$-submodules. We denote by $V_o$
the submodule on which the action is trivial; while $V^k$ are the
non-trivial submodules of dimension one and $V_k$ the ones that
have dimension bigger than $1$. Using Theorem 1.3, we have that
$\dim V_o \leq 2$ and $r+j \leq 2$. We distinguish some cases
according to the
dimension of $V_o.$ \\
When $\dim V_o=n-k=2$ we have $r=0$ and $j=0$, recalling that
$2k\leq n,$  we conclude that $k=2$. In this case $\mf l=\su (2)$
because $\li 2=0;$ since the scalars cannot be removed when $\su
(2)$ acts diagonally, the action is not
multiplicity free.\\
If $\dim V_o=1$ either $r=1$ and $j=0$ or $r=0$ and  $j\leq 1.$
\\Since the trace $tr(\rho_2(X))=0$ for all $X\in \mf{l},$ if
$r=0$ then $j=0$  and in particular $k=1,$ ( $2k\leq n=k+1$)
contrary to our assumption. We have the following possibilities,
from Tables IIa and IIb.
$$
\begin{array}{|c|c|c|c|}
n. & \mf{l}_1  &  \mf{l}_2     & {\rm conditions} \\ \hline 1.  &
\su (k)   & \su (n-k-1)      &  n \geq k +3       \\ 2.  & \su(2)
& \mf{sp} (\frac{n-3}{2}) & n \geq 7, k=2    \\
\end{array}
$$
From this we see that $\mf l\subseteq \mf z+ \li 1+\li 2.$ Note
that when $n>2k+1$ in case $1$ the scalars can be reduced while in
case $2$, scalars can never be reduced. Summing up we conclude
that only  $\mf z+ \su (k) + \su (n-k-1)$   on $ Gr (k,n),$
$ n >2k+1,$ gives rise to a coisotropic action.\\
The last case arises when $\dim V_o=0.$ We have $r+j=2$, otherwise
we fall again in the  case in which both the actions were
irreducible.\\
If $r=2$ by Tables IIa and IIb we have
$$
\begin{array}{|c|c|c|c|}
n. & \li 1 & \li 2 & {\rm conditions} \\ \hline
1.  & \su(2) & \su (p) + \su (q) & p+q=n-2,\ \ p,q \geq 2 \\
2. & \su(2) & \su (p) + \asp (q) & p+2q=n-2,\ \ p \geq 2,\ q > 1 \\
3. & \su (2) & \asp (p) + \asp (q) & 2p+2q=n-2, \ \ p,q > 1 \\
\end{array}
$$
Note that the centralizer is $\mf {z} + \R$, where $\R$ denotes
the centralizer of $\li 2$ in $\su (n-2).$ We recall briefly the
action of $\mf{z}$ and $\R$ on the slice that is given by $$
(\pi^* \otimes V_1) \oplus  ( \pi^* \otimes V_2). $$ $\mf{z}$ acts
as $(e^{i ( \frac{\phi}{2} + \frac{\phi}{n-2})}, e^{i (
\frac{\phi}{2} + \frac{\phi}{n-2}) } )$ while $\R$ acts as $
(e^{-i \frac{\psi}{p} }, e^{i \frac{\psi}{q} }).$\\ By  Table IIa
the scalars, in the first case, can be removed if $p$ and $q$ are
bigger than $2$ while they cannot be reduced if $p=q=2.$ If $p=2$
and $q>2$  the scalars cannot be reduced if the action on the
first factor is  trivial: this happens when we consider the line $
\psi = \frac{n}{n-2} \phi,$ which one can show  easily that
corresponds to the centralizer of $\su(n-4)+\su (4)$ in $\su
(n).$\\ In the second case the scalars cannot be reduced if $p=2.$
If $p
>2$ then the scalars are needed  if the action on the second
factor is trivial, that corresponds to the line $\psi=
\frac{-nq}{2(n-2)} \phi.$\\ Finally in the third case the scalars
can never be reduced.\\ Here we give the corresponding minimal
coisotropic actions. In what follows we will denote by $\R
({\alpha} )$ the lines different from $y= \alpha x$ where $(x,y)
\in \mf{z} + \R.$ $$
\begin{array}{|c|c|c|c|}
n. & \h  &  M   & {\rm conditions} \\ \hline
1.  & \mf z+ \su (2) + \R + \su (2) + \su(2)   & Gr (2,6)  & \\
2.  & \su (2) + \su (p) + \su(q)   & Gr (2,n)  &  p \ {\rm and} \ q >2,
p+q=n-2 \\
3.  & \R({\frac{n}{n-2}}) + \su (2) + \su (2) + \su(n-4)   & Gr
(2,n)  &  n>6
\\
4.  & \R({\frac{-nq}{2(n-2)}}) + \su (2) + \su (p) + \asp (q)
& Gr (2,n) & p>2,\ q>1, \ n=2+p+2q \\
5.  & \mf z + \su (2) + \R + \su (2) + \asp (n-2)
& Gr (2,2n)) & n >3 \\
6.  & \mf z  + \su(2) +\R + \asp (p) + \asp (q)  & Gr (2,2n) & p>1,\ q>1, \
2p+2q=n-2\\
\end{array}
$$
Suppose now $j=1$ and $r=1$. We obtain
$$
\begin{array}{|c|c|c|c|}
n. & \mf{l}_1  &  \mf{l}_2     & {\rm conditions} \\ \hline 1.  &
\su (k)   & \su (n-k-1)+ \R      &  n \geq k +3  \\ 2.  & \su(2) &
\mf{sp} (\frac{n-3}{2})+\R & n \geq 7, k=2    \\
\end{array}
$$
Note that the centralizer is again $\mf {z} + \R$, where $\R$ is
the centralizer of $\li 2$ in $\su (n-k).$ Clearly $\mf l$ can be
$ \mf z + \R + \su (k)^{\Delta} $ but this case does not hold for
dimensional reasons. The scalars acts as follows: $\mf{z}$ acts as
$(e^{i ( \frac{\phi}{k} + \frac{\phi}{n-k})}, e^{i (
\frac{\phi}{k} + \frac{\phi}{n-k}) } ),$ and $\R$ as  $ (e^{-i
\psi }, e^{i \frac{\psi}{n-k-1}}).$ Hence, with the same arguments
used above, we get
$$
\begin{array}{|c|c|c|c|}
n. & \h  &  M   & {\rm conditions} \\ \hline 1.  & \mf z + \su (k)
+ \R + \su (k-1)   & Gr (k,2k)  & k\geq 3   \\ 2.  & \mf z + \su
(k) + \R + \su (k)   & Gr (k,2k+1)  & k\geq 2   \\ 3. &
\R({\frac{n}{(n-k)k}}) + \su(k) + \su (n-k-1)  & Gr (k,n) & n >2k
+1
\\ 5.  & \mf z + \R+\su(2) + \mf{sp}(n-1)   & Gr(2,2n+1) & n >
2 \\
\end{array}
$$ Finally if $j=2$, we get that $n-k=2$, hence $k=2,$ $\li
2=\mf{t}_1\subseteq \su (2)$ is a torus and we have the
corresponding diagonal action of $\mf{su} (2)$ on the slice, given
by $\C^2 \otimes \C^2.$ Note that, in this case, the centralizer
coincides with $\mf {z} + \mf {t}_1$, hence the action of
$\mf z +\mf {t}_1+\su (2) $ on $Gr(2,4)$ is
multiplicity free but the scalar cannot be reduced.

Now assume $\li 1$ acts reducibly and $\li 2$ acts irreducibly.
Following the same procedure used in the symmetric case we find
either the dual cases, i.e. where one interchange the roles of $k$
and $n-k$, or cases already found. Only when, with the same
notations used above, $ ( \C^k )^* = W_o \oplus W^1 \oplus \cdots
\oplus W^m \oplus W_1 \oplus \cdots \oplus W_s$ and $\dim W_o=0,$
$m=2$ and $s=0$, we get another different coisotropic action.
Clearly, again by Tables IIa and IIb, $k=2$ and on the slice we
have $\su (n-2)$ or $ \asp (\frac{n-2}{2})$ acting diagonally on
$\C^{n-2}\otimes \C^{n-2}.$ The scalars can be reduced if $n \neq
4$ in the first case, while in the second case the scalars are
needed. We note also that $\mf{l}_1$ is a torus, denoted by $\mf
t_1,$  in $\su (2).$ Hence we have, with the same notation used
above, the following multiplicity free actions.
$$
\begin{array}{|c|c|c|c|}
n. & \h & M & {\rm conditions} \\ \hline
1.  & \mf {t}_1 + \su(n-2)
& Gr (2,n) & n > 4 \\
2. & \mf z + \mf {t}_1 + \asp (n-1) & Gr
(2,2n) & n \geq 2
\\
\end{array}
$$
\\
{\bf (3)} $\li 1$ and $\li 2$ act both reducibly on $(\C^k )^*$
and $\C^{n-k}$ respectively. We can split $ (\C^k)^* = W_o \oplus
W^1 \oplus \cdots \oplus W^m \oplus W_1 \oplus \cdots \oplus W_s $
as $\li 1$-irreducible submodules and $ \C^{n-k} = V_o \oplus V^1
\oplus \cdots \oplus V^j \oplus V_1 \oplus \cdots \oplus V_r $ as
$\li 2$-irreducible submodules. First of all note that $r \neq 0$
or $s \neq 0,$ otherwise the dimensional condition is not
satisfied. We here investigate the possible dimensions of $V_o$
and $W_o$. Assume that $V_o, W_o \neq 0;$ then $\li 1 + \li 2$
acts trivially on $V_o \otimes W_o$ and $\mf l$ acts
coisotropically, so that $\mf l$ has no trivial projection along
$\mf z$ which acts on $V_o \otimes W_o.$ Since $\dim \mf z=1,$ we
must have $\dim (V_o \otimes W_o ) \leq 1.$ On the other hand
$\dim W_o,\dim V_o \leq 2.$ Indeed, assume by contradiction, $\dim
W_o
>2.$ Clearly $r=0,$ otherwise we have too many terms in the slice;
on the other hand by dimensional reasons $s$ has to be greater or
equal than $1$. Recall also that $tr( \rho_2 (X) )=0, \forall X
\in \h$ so $\dim V_o=0$ and $j=2$ or $\dim V_o
>0$ and $j=0$. The first case is excluded otherwise $k=2$ that is
impossible since $k\geq \dim W_o>2$. Hence $\dim V_o \neq 0$ so
$\dim (W_o \otimes V_o) \geq 2$ that is a contradiction. \\ Note
also that we have a symmetry between $W_o$ and $V_o$ so we have
only four cases. Our aim is to prove none of these gives a
coisotropic action. We will analyze the one in which $\dim V_o =
\dim W_o=1$ because the others are quite similar. We know that $r$
or $s$ must be positive. Assume, for example, $r \geq 1.$ Since
$\rho_1 (l) \subseteq \su (k)$ then, again by the condition on the
trace, $m=0.$ In particular $s$ has to  be greater, as $k> 1,$ or equal to $1$ so
$j=0$ as before. Now in the slice representation appears  the
indecomposable factor $(W_1 \otimes V_o) \oplus ( V_1 \otimes W_o)
\oplus ( W_1 \otimes V_1),$ hence the action cannot be
multiplicity free. Clearly the case $s \geq 1$ is exactly the
same.
%
%
%
%
%
%
%
%%%%%%%%%%%%%%%%%%%%%%%%%%%%%%%%%%%%%%%%%%%
%
% IL CASO RIDUCIBILE IN GENERALE SE k<l e k<n-l
%
%%%%%%%%%%%%%%%%%%%%%%%%%%%%%%%%%%%%%%%%%%%%%%%%%%%%%%%%%%%%%%%%%%%%%%%%%%%%%%%%%%%%%%%%%%%%%%%%
\subsubsection{The case $k<l,\ k<n-l$}
Let $L \subset \SUr{l}{n-l}$. Let $\pi_1$ be a complex $k$-plane
in $\C^l$ and let $\pi_2$ be a complex $k$-plane in $\C^{n-l}$.
Clearly $K \pi_i$ are complex orbits so if $L \subseteq K$ acts on
$Gr (k,n)$ coisotropically then it acts also on $K \pi_i$,
$i=1,2$, coisotropically (see \cite{Hw} Restriction Lemma). \\
With the same notations used in the previous case we denote by
$\li 1$ (resp. $\li 2$) the projection of $\mf l$ on $\su (l)$
(resp. $\su (n-l)$).  We have the following possibilities:
\begin{description}
\item[(1)] $\mathfrak{l}_i,\ i=1,2$ acts transitively on $K\pi_i$ ;
\item [(2)] $\mathfrak{l}_i$ acts transitively and $\mf{l}_j (\mathfrak{l})$
acts coisotropically ($ i \neq j \ mod\, 2$);
\item [(3)] $\mf{l}_i,\ i=1,2$ acts coistropically.
\end{description}
{\bf (1)} $\mathfrak{l}_1= \mathfrak{su}(l)$ and $
\mathfrak{l}_2=\mathfrak{su}(n-l).$ The simply connected group
with Lie algebra $\mathfrak{l}=\li 1+\li 2$ is $L=\SU(l) \times
\SU(n-l)$ and $L \pi_1$ is a complex orbit. Furthermore, it is
easy to see that the slice representation is given by $S= \pi_1^*
\otimes \C^{n-l}$. On the slice $S$ we have, looking at Tables Ia,
Ib, the irreducible multiplicity free action of $\su (k)+ \su
(n-l)$, where $k\neq n-l,$ and $k,n-l\geq 2.$ Hence, under these
conditions, we get that this action is multiplicity free. Note
that when $l=n-l$ then  $\h$  should be also
$\R+{\mf{su}(l)}^\Delta$, but this case is not coisotropic for
dimensional reasons.
\\
\\
\n {\bf (2.a)} $ \mathfrak{l}_1 = \mathfrak{su}(l)$ and $
\mathfrak{l}_2=\mathfrak{so}(n-l)$ (or $\mathfrak{sp}(n-l)$). We
can take a complex $\mf{so}(n-l)$ (resp. $\asp ({n-l})$)-orbit in
$\C^{n-l}$ and it is easy to see that the slice representation
contains the  submodule $S^2(\pi_2 ) \oplus \pi_2^* \otimes \C^l$
(resp. $\Lambda^2(\pi_2 ) \oplus \pi_2^*  \otimes \C^l$). The
stabilizer's factor $\U(k)$ acts  diagonally on  $S^2(\pi_2)$
(resp. on $\Lambda^2 (\pi_2)$ ) and on $\pi_2^*$. This kind of
indecomposable action does not appear in Table II, hence these
actions fail to be coisotropic.\\
\\
\n {\bf (2.b)} $\mathfrak{l}_1 = \mathfrak{su}(l)$ and
$\mathfrak{l}_2 \subseteq \R + \su (p)+ \su (q),$ where $ p+q=n-l$
and $\R$ is the centralizer of $\li 2$ in $\su (n-l).$ The slice
representation relative to the complex orbit $\SU(l) \pi_1$ is
given by $S=\pi_1^* \otimes \C^{n-l}$. We can decompose $\C^{n-l}$
as $ \mathfrak{l}_2$-irreducible submodules as follows
$V_o \oplus V^1 \oplus \cdots \oplus V^j
\oplus V_1 \oplus \cdots \oplus V_r,
$
with the same notations used in the fixed point case. Note that
again $\dim V_o \leq 2$ and $r+j \leq 2.$ We will follow exactly
the same procedure used when $k=l$ distinguishing three different
cases, namely when $\dim V_o=2,1$ or $0.$
\\
When $\dim V_o=2$ we
have $\C^{n-l} =V_o$ that implies, since $k\neq 1,$ $k=2=n-l,$
which is not our case, since $k<n-l.$
\\
If $\dim V_o=1$ it is easy to see that $\C^{n-l}$ splits as the
direct sum $\C \oplus V_1$, i.e $j=0$ and $r=1.$ The slice becomes
$S={\pi_1}^*\oplus {\pi_1}^*\otimes V_1,$ hence on the slice $\R_s
+\su(k)+\su(n-l-1)$ or $\R_s +\su(2)+\mf{sp}((n-l-1)/2)$ act,
where $\su(k)+\R_s$ acts diagonally on ${\pi_1}^*\oplus {\pi_1}^*$
and where $\R_s \subseteq \su (l)$ is the determinant condition of
the stabilizer's factor $\SUr{k}{l-k}$ since $k<l.$ The
centralizer of $\mf{l}$ coincides with $\mf c+\R_s$, where the
centralizer $\mf c$ of $\su (l) + \su(n-l)$ in $\su (n)$  acts as
$(e^{i(\frac{n \phi}{(n-l)l})},e^{i(\frac{n \phi}{(n-l)l})})$,
while $\R_s$ acts as $(e^{i(\frac{ \psi}{k})},e^{i(\frac{
\psi}{k})}):$ that corresponds to a single scalar action. Note
that, only in case $1,$ when $n>k+l+1,$ one can reduce the
scalars. Hence we find only the coisotropic action of  $\su(l) +
\su (n-l-1)$ on $Gr (k,n),$ where $ n >k+l +1.$ \\ Suppose now
$\dim V_o=0.$ Assume first $r=2.$ The slice $S$ splits as
${\pi_1}^*\otimes V_1\oplus {\pi_1}^*\otimes V_2,$ hence $k=2$ and
on the slice we have the following multiplicity free actions $$
\begin{array}{|c|c|c|c|}
n. & \mf{h}_1 & \li 2 & {\rm conditions} \\ \hline
1.  & \R_s+ \su(2)& \su (p) + \su (q) & p+q=n-l,\ \ p,q \geq 2 \\
2. & \R_s+ \su(2) & \su (p) + \asp (q) & p+2q=n-l,\ \ p \geq 2,\ q
\geq 1 \\
3. & \R_s+ \su (2) & \asp (p) + \asp (q)& 2p+2q=n-l, \ \ p,q \geq
1
\\
\end{array}
$$
where $\mf{h}_1$ denotes the factor of the stabilizer of $\pi_1$
that acts on ${\pi_1}^*.$ Note that the centralizer of $\mf{l}$ in
$\mf{k}$ is $\mf {c} + \R_s+\R .$ Here we give the corresponding
minimal coisotropic actions,
$$
\begin{array}{|c|c|c|c|}
n. & \h  &  M   & {\rm conditions} \\ \hline 1.  &  \su (l) + \R +
\su (2) + \su(2)   & Gr (2,n)  &    n=l+4 \\ 2.  & \su (l) + \su
(p) + \su(q)   & Gr (2,n)  &  p,q >2,  p+q=n-l \\
3.&\su (l) + \su (2) + \su(n-l-2)   & Gr (2,n)& k<l,n-l \\
4. & \su(l) +\su (p) + \asp (q)  & Gr (2,n) & p>2, p+2q=n-l \\
5.& \su(l)+\R+\su (2)+\asp (q)&Gr(2,n) & 4\leq 2q=n-l-2\\
6.  & \su(l) +\R + \asp (p) + \asp (q)  & Gr (2,n) & 2p+2q=n-l,\; p,q\geq 2 \\
\end{array}
$$
Suppose now $j=1$ and $r=1$. We obtain, with the
same arguments  used in the fixed point case, and recalling
that the centralizer is again $\mf {c} + \R$
the following minimal coisotropic actions:$$
\begin{array}{|c|c|c|c|}
n. & \h  &  M   & {\rm conditions} \\ \hline 1.  &  \su (l) + \R +
\su (n-l-1)   & Gr (k,n)  &  k=n-l-1 \geq 2\\ 2. &  \su(l) + \su
(n-l-1)  & Gr (k,n) & n >k +l+1 \\ 3.  & \R+\su(l) + \mf{sp}(n-1)
& Gr(2,2n+l-1) & n > 2 \\
\end{array}
$$
Finally if $j=2$ and $r=0$   we get that $n-l=2$, hence $k=1$,
that is not our case. \\
\\
\n
{\bf (3)} Note that, by {\bf (2.a)}, we can restrict our analysis
to the case in which the two projections are contained into
reducible subgroups; that is $ \mathfrak{l}_1 \subseteq \R
+\mathfrak{su}(p_1) + \mathfrak{su}(q_1),\ p_1 + q_1=l$ and
$\mathfrak{l}_2 \subseteq \R + \mathfrak{su}(p_2) +
\mathfrak{su}(q_2),\ p_2 + q_2=n-l.$ We want to prove that this
action fails to be multiplicity free.\\ Take $\pi_1$ a complex
$k_1$-plane in $\C^{p_1}$ and $\pi_2$ a complex $k_2$-plane in
$\C^{q_1}$, where $k_1+k_2=k.$ The orbit $\Su (\U (p_1) \times \U
(q_1)) (\pi_1 \oplus \pi_2 )$ is complex and it is easy to check
that the slice is given by
$$ (\pi_1^*
\otimes \pi_2^{\perp} )\oplus (\pi_1^* \otimes \C^{p_2}) \oplus
(\pi_1^* \otimes \C^{q_2})  \oplus (\pi_2^* \otimes
\pi_1^{\perp})\oplus (\pi_2^* \otimes \C^{q_2})\oplus (\pi_2^*
\otimes \C^{p_2}).
$$
where $\pi_1^{\perp}$ (resp. $\pi_2^\perp)$ denotes the orthogonal
complement of $\pi_1$ (resp. $\pi_2$) in $\C^p$ (resp. $\C^q$). \\
Note that only the case $k_1=k_2=1$ is admissible since, otherwise
there would be at least three terms on which for example
$\mathfrak{su}(k_1)$ (if $k_1>1$) acts. Then, under this
assumption, there must be, for example $p_1>1,$ otherwise
$p_1+p_2+q_1+q_2=n=4,$ and the corresponding action is not
coisotropic for dimensional reasons. Take a $2$-plane $\pi$ in
$\C^{p_1}$, then in the slice appears the term $ \pi^* \otimes
\C^{p_2}\oplus \pi^* \otimes \C^{q_1} \oplus \pi^* \otimes
\C^{q_2},$ hence the action is not coisotropic.
%%%%%%%%%%%%%%%%%%%%%%%%%%%%%%%%%%%%%%%%%%%%%%%%%%%%%%%%%%%%%%%%%%%%%%%%%%%%%%%%%%%%%%%%%%%%%%%%%%%%%%%%%%%%%%%%%%%%%%%
%
% IL CASO  RIDUCIBILE  con k>n-l
%%%%%%%%%%%%%%%%%%%%%%%%%%%%%%%%%%%%%%%%%%%%%%%%%%%%%%%%%%%%%%%%%%%%%%%%%%%%%%%%%%%%%%%%%%%%%%%%%%%%%%%%%%%%%%%%%%%%%%
%
% SE LA PROIeZIONE AGISCE TRANSITIVAMENTE
%%%%%%%%%%%%%%%%%%%%%%%%%%%%%%%%%%%%%%%%%%%%%%%%%%%%%%%%%%%%%%%%%%%%%%%%%%%%%%%%%%%%%%%%%%%%%%%%%%%%%%%%%%%%%%%%%%%%%%
\subsubsection{The case $k<l, \ k>n-l$}
We can take a $k$-plane $\pi\subseteq \C^l$ and determine the
complex $K$-orbit. By the restriction lemma we have that the
projection $\rho_1(\mf{l})=\mf{l}_1$ must act coisotropically on
the complex orbit $Gr (k,l).$ Therefore we may have two possible
situations, namely when the action is transitive, i.e
$\mf{l}_1=\mf{su}(l)$, or when $\mf{l}_1$ acts coisotropically on
$Gr (k,l)$.\\
\\
{\bf (1)} Suppose $\mf{l}_1=\mf{su}(l).$ Note that $\mf{l}_1$ acts
on the first factor of  the slice $S=\pi^*\otimes \C^{n-l}$
irreducibly; while the second factor can split as the direct sum
of irreducible $\mf{l}_2$-submodules.
Following the same procedure we have used when both $l,n-l> k,$
and with the same notations, we have $\dim V_o\leq 2,$ and
$\C^{n-l}$ can split as at most two summands.\\ When $\dim V_o=2$
the action is not coisotropic, as the diagonal action of
$\mf{su}(k)$ needs  the scalars indeed the center acts as
$z\rightarrow(z^a,z^a).$ \\ If $\dim V_o=1$ then $r=1$ and $j=0$,
or $j=r=0;$ in the first case the only possibility  is the action
of $\mf{su}(k)+\mf{su}(n-l-1)$, $n-l \geq 3$, on the slice; note
that the scalars are removable if and only if $n-l-1>k$ which is
not our case, hence the action is not multiplicity free. In the
second case, by Table Ia, we get the multiplicity free action of
$\su (n-1).$ \\ Finally if $\dim V_o=0$ we obtain the same cases
found in {\bf (2.b)} (when $\dim V_o=0$). However, note that the
groups arising when  $r=2$ and $j=0$ have to be excluded since
otherwise, by Tables IIa, IIb, $k=2$ and $2>n-l=\dim V_1+\dim
V_2\geq 4$ that gives rise to a contradiction. On the other hand
the case arising when $j=2$ and $r=0$ has to be considered. In
this situation  $\mf{l}_2=\mf{t}_1$ and $l=n-2$.\\ Summing up,
when the $\li 2$-action is reducible, we get the following
coisotropic actions :
$$
\begin{array}{|c|c|c|c|}
n. & \h  &  M   & {\rm conditions} \\ \hline 1.  &  \su (l) + \R +
\su (n-l-1)   & Gr (k,n)  &  n \geq l +5, k>n-l  \\ 2.&\mf
{t}_1+\su (n-2)& Gr(k,n)& k>2\\
\end{array}
$$
When the $\li 2$-action is irreducible, with the same arguments
that we have used in the fixed point case, we have the following
multiplicity free actions.
$$
\begin{array}{|c|c|c|c|}
n& \mf{l} & M & {\rm conditions}  \\ \hline
1.&\su (n-1)    & Gr(k,n) &  n \geq 4 \\
2.& \mf{su}(l)+\mf{su}(n-l) & Gr(k,n)  &  k, n-l\geq 2 \\
3.&\mf{sp}(2)+\mf{su}(n-4) & Gr (k,n)  & 5 < k\leq \frac{n}{2} \\
\end{array}
$$
%%%%%%%%%%%%%%%%%%%%%%%%%%%%%%%%%%%%%%%%%%%%%%%%%%%%%%%%%%%%%%%%%%%%%%%%%%%%%%%%%%%%%%%%%%%%%%%%%%%%%%%%%%%%%%%%%%%%%%%%%%%
%
%
%             SE NON AGISCE TRANSITIVAMENTE 1 MODO
%
%
%%%%%%%%%%%%%%%%%%%%%%%%%%%%%%%%%%%%%%%%%%%%%%%%%%%%%%%%%%%%%%%%%%%%%%%%%%%%%%%%%%%%%%%%%%%%%%%%%%%%%%%%%%%%%%%%%%%%%%%%%%%
\\
\n {\bf (2)} Suppose now that $\mf{l}_1$ acts coisotropically on
$Gr(k,l).$ We will test all the maximal  subgroups  of $\SU(l)$ that act
multiplicity free on $Gr (k,l)$ and then, according to the analysis that we have already done, we will go through all their subgroups still acting coisotropically. The idea will be always the same,
we choose a complex orbit, we find the slice representation and we
determine under which conditions the action appears in Tables Ia
or IIa, IIb.
\begin{description}
\item[]$\mf{l}_1=\mf{so}(l).$ The slice
$S=S^2(\pi^*)\oplus (\pi^*\otimes \C^{n-l})$ does not appear in
Tables IIa, IIb; hence the action is not coisotropic.
\item[]$\mf{l}_1=\mf{sp}(l/2).$ The slice is
$S=\Lambda^2(\pi^*)\oplus (\pi^*\otimes \C^{n-l})$ where $\U(k)$
acts diagonally on $\Lambda^2(\pi^*)\oplus \pi^*$. It appears in
Table IIa or IIb iff $n-l=1$, i.e. $l=n-1.$ The scalars act on the
slice as $z\mapsto (z^a,z^a)$. Therefore, by Tables IIa and IIb,
$\mf{sp}((n-1)/2)\subseteq \mf{su}(n-1)$ acts coisotropically on
$Gr (k,n)$ for $k\neq 3$;  if $k=3$,  the slice becomes
$S=\pi^*\oplus \pi^*$ and $\R+\asp (n)$ act coisotropically on
$Gr(3,2n+1).$
 Note that no
subgroup of $\Sp(n-1/2)$ can act coisotropically on the complex
Grassmannian.
\item[]$\mf{l}_1=\R+\mf{su}(p)+\mf{su}(q)$ where $p+q=l$ and $\R$
is the centralizer of $\mf{su}(p)+\mf{su}(q)$ in $\mf{su}(l)$. We
take $k_1, \ k_2$ such that $k_1 \leq p, \ k_2 \leq q$ and $k=k_1
+ k_2$. Let $\pi_1$ and $\pi_2$ be $k_1$- and $k_2$-planes in
$\C^{p}$ and $\C^{q}$ respectively; let $\pi= \pi_1 \oplus \pi_2$
be the $k$-complex plane in $\C^l.$ The slice is
$$ S=({\pi_1}^*\otimes \C^{n-l})\oplus ({\pi_2}^*\otimes
\C^{n-l})\oplus ({\pi_1}^*\otimes \pi_2^{\perp})\oplus
({\pi_2}^*\otimes \pi_1^{\perp}). $$
We are going to prove that, for all $k$, $n-l$ has to be $1$ and
$\li 2 =0.$ Indeed, if $k=2$ then, since  $k>n-l$ then $n-l=1$.
For $ k \geq 3$ then the group acts on the slice coisotropically
only if $n-l=1$, i.e. $p+q=n-1,$ note that if it is not the case,
there would be too many terms.\\ 
By a straightforward calculation
we can always take $p$ or $q$ bigger than $k,$ otherwise we fall
again in the fixed point case. Assume, for example, $p>k$ and take
a complex $k$-plane $\pi \subset \C^p$ and consider the complex
orbit through $\pi.$ The slice representation is given by
$ (\pi^* \otimes \C^q) \oplus (\pi^* \otimes \C)$ where $\R_s+ \su
(k)$ acts diagonally on $\pi^*$; here $\R_s$ denotes the
centralizer of $\su (k) + \su (p-k)$ in $\su (p)$. \\ The
centralizer,  $\mf c ,$ of $\su (n-1)$ in $\su (n)$, acts on the
slice as $(1,e^{-i\frac{n\phi}{n-1}} );$
 while $\R$  as $ (e^{-i
\frac{ \psi}{p} } e^{-i \frac{\psi}{q} }, e^{-i \frac{ \psi}{p}})$
and $\R_s$ as $ (e^{i \frac{\theta}{k} }, e^{i \frac{ \theta}{k}
}).$ \\ If $q=1$ we shall assume $k \geq 3$ since if it is not we
fall again in the fixed point case. The slice is given by $\pi^*
\oplus \pi^*$ and the scalars act on the slice
$ (e^{-i \frac{ \psi}{p}} e^{-i \psi} e^{i \frac{ \theta}{k}},
e^{-i \frac{ n \phi}{n-1}} e^{-i \frac{ \psi}{p} } e^{i \frac{
\theta}{k}} ).$
We know that this action fails to be multiplicity free when $a=b.$
We conclude, when $q=1,$  that the action is coisotropic if and
only if $\psi\neq \frac{n \phi}{ n-1}.$ Note also that the line
$\psi= \frac{n \psi}{ n-1}$ corresponds to the centralizer of
$\su(2) + \su (n-2)$ in $\su (n).$ In particular, if $\mf {t}_1$
denotes  a torus in $\su (2),$ then the action is multiplicity
free if and only if the projection of the scalar has a no trivial
component along $\mf t_1. $ \\ Now, suppose that $q >1.$ The slice
becomes $(\pi^* \otimes \C^q) \oplus (\pi^* \otimes \C) $ and the
action of the scalars on the slice is $ (e^{-i \frac{ \psi}{p} }
e^{-i \frac{\psi}{q} } e^{i \frac{\theta}{k} }, e^{-i \frac{ n
\phi}{n-1} } e^{-i \frac{ \psi}{p} } e^{i \frac{ \theta}{k} }).$
\\ By Tables IIa, IIb we  can distinguish some cases according to
the relations between $k$ and $q.$ \\ Firstly if $q>k$ one can
reduce the scalars if and only if the action on the second factor
is non trivial and this holds since $\R_s$ acts. \\ When $q=k$ we
fall again in the fixed point case. \\ In  case $q=k-1$ the
scalars cannot be reduced and it is easy to check that the image
of the scalars's action is one dimensional if and only if $\psi=
\frac{(k-1)n \phi}{n-1}.$ \\
%that correspond to the centralizer of $\su (k+1) + \su (p)$ in $\su (n).$
%Indeed on $\C^q$ that line acts as $e^{-i \phi}.$
Finally if $q < k-1$ the action fails to be multiplicity free if
and only if the image of the scalars is one dimensional and on
both factors the action is the same. This happens on the line
$\psi = \frac{qn \phi}{n-1}.$
%that corresponds to the centralizer
%of $\su (q+1) + \su (n-q-1)$ in $\su (n).$
Summing up we have the following minimal coisotropic actions
$$
\begin{array}{|c|c|c|c|}
n.& \h & M & {\rm conditions} \\ \hline
1.& \mf t_1 + \su (n-2) & Gr (k,n) &  k>2 \\
2.& \su (p) + \su (q)  & Gr (k,n)  &  p+q=n-1,\  p,q>k \\
3.& \R ( \frac{ (k-1) n }{n-1} ) + \su (k-1) + \su (n-k) & Gr(k,n) & k \geq 3 \\
4.& \R (\frac{qn}{n-1}) + \su (p) + \su (q)& Gr (k,n) & p+q=n-1,\
p>k,\ k \geq q+2. \\
\end{array}
$$
\end{description}
We shall investigate if we can go under $\R+\su(p)+\su(q)$.\\ We first consider the maximal subalgebras of $\su(q)$. Choose a $k$-plane $\pi$ in $\C^p$ as above. The slice, corresponding to the orbit through $\pi$  is given by $S=\pi^*\otimes \C^q \oplus \pi^*\otimes \C$. By Tables IIa and IIb we argue that no subalgebras give rise to a coisotropic action.\\
Let now consider the maximal subalgebras of $\su (p)$.\\ We take the subalgebra $\mf{so} (p)$ (resp. $\asp (p/2)$). Take the $k$-plane $\pi$ such that the orbit through it is complex. The slice becomes $S=S^2(\pi^*)\oplus (\pi^*\otimes \C)\oplus (\pi^*\otimes \C^q)$ (resp. $S=\Lambda^2(\pi^*)\oplus (\pi^*\otimes \C)\oplus (\pi^*\otimes \C^q)$) and the corresponding actions on it do not appear in Tables IIa and IIb. The case of a simple subalgebra can be treated
likewise.\\ Let now investigate $\R+\su(p_1)+\su(p_2)\subset \su(p)$.
We
take $k_1, \ k_2$ such that $k_1 \leq p_1, \ k_2 \leq p_2$ and $k=k_1
+ k_2$. Let $\pi_1$ and $\pi_2$ be $k_1$- and $k_2$-planes in
$\C^{p_1}$ and $\C^{p_2}$ respectively; let $\pi= \pi_1 \oplus \pi_2$
be the $k$-complex plane in $\C^p.$ The slice  contains the submodules
$({\pi_i}^* \otimes \C)\oplus ({\pi_i}^*\otimes
\C^q)\oplus ({\pi_i^*} \otimes \pi_j^{\perp}) $ with $i\neq j$, hence the corresponding action is not coisotropic.\\
Finally the case $\mf {l}_1=\R+\R+\su(p)^\Delta$ can be excluded for dimensional reasons.
%%%%%%%%%%%%%%%%%%%%%%%%%%%%%%%%%%%%%%%%%%%%%%%%%%%%%%%%%%%%%%%%%%%%%%%%%%%%%%
%
%        POLARITA'
%
%
%%%%%%%%%%%%%%%%%%%%%%%%%%%%%%%%%%%%%%%%%%%%%%%%%%%%%%%%%%%%%%%%%%%%%%%%%%%%%%
\section{Polar actions on Complex Grassmannians}
%inizio
In this section we study which coisotropic actions
 are polar. Firstly observe that the reducible actions arising from Table IIa
and IIb are not polar; this can be easily deduced as an application of
Theorem $2$ (pag. 313) \cite{Bergmann}; while, see \cite{He} and \cite{Kol},
in the irreducible case we know  that $\mf{so}(n)$
on $Gr(k,n),$ $\asp (n)$ on
$Gr(k,2n)$ and  $\su (l)+\su (n-l)$,  for $l\neq n-l$ on $Gr(k,n)$
 give rise to polar actions.
The last two cases $\mf{spin}(7)$ on $Gr(2,8)$ and $\asp (n)$ on
$Gr(k,2n+1)$ can be excluded. The first one fails to be polar
since $\mf {spin}(7)$ has the same orbits of $\mf{u}(3)$ which is
not polar for \cite{He}.\\ Let now consider the action of $\asp
(n)$ on $Gr (k,2n+1).$ For $k \geq 3,$ using Theorem $2$
\cite{Bergmann}, one can easily prove that the action is not
polar. Note that if $k=4$ the action on $S= \Lambda^2 (\pi^*)
\oplus \pi^*$ appears in \cite{Bergmann} but as real
representation. \\Now we consider the case $k=2.$ As $\asp
(n)$-module  $\C^{2n+1}= \C \oplus \C^{2n}$ and take $\pi$ in
$\C^{2n}$ such that the orbit through $\pi$ is complex. The normal
space $N=\Lambda^2(\pi^*) \oplus \C^2$ is acted on by $\Sp
(n)_{\pi}=\U(2)\times \Sp (n-2)$ so that the action  has
cohomogeneity $2$. If the $\Sp (n)$ action were polar, using
Proposition 1 \cite{Bergmann} a section in $N$ can be taken as the
direct sum of a section for the action of $\U (2)$ on $\Lambda^2
(\pi^*)$ plus a section of $\U (2)$ on $\C^2.$ With this remark
any subspace $\mf m \subseteq N$ generated by $v_1 \in \Lambda^2 (
\pi^*)$ and $v_2$ in $T_{\pi} Gr (2,2n)^{\perp} \subseteq T_{\pi}
Gr (2, 2n+1)$ would be the tangent space to a totally geodesic
submanifold. A direct inspection taking, $$ v_1=\left(
\begin{array}{cccccc}
0&0&0&1&\ldots&0  \\
0&0&-1&0&\ldots&\vdots \\
0&1&0&0&\ldots&\vdots\\
-1&0&0&0&\ldots&\vdots \\
\ldots&\ldots&\ldots&\ldots&\ldots&\vdots\\
0&\ldots&\ldots&\ldots&\ldots&0 \\
\end{array}
\right) v_2=\left(
\begin{array}{ccccc|c}
 0 & 0&\ldots & 0&\ldots&1     \\
 0   & 0 & \ldots &0&\ldots& 0  \\
 \ldots & \ldots & \ldots &0& \ldots&\vdots   \\
 0&0&0&0&\ldots&\vdots\\
 \ldots&\ldots&\ldots&\ldots&\ldots&0 \\
\hline -1   & 0& \ldots &\ldots&0& 0     \\
\end{array}
\right) $$ shows, together with Theorem $7.2$ pag. 224 \cite{Hel}
on Lie triple systems, that the section  $\Sigma= \exp (\mf m)$ is
not totally geodesic, hence the action cannot be polar.\\ Finally
observe that the groups $K$ that are listed in Table $3$ act also
hyperpolarly on complex Grassmannians, since  $K$ and
$\SUr{l}{n-l}$ are symmetric subgroups of $\SU (n)$.
\newpage
\section{Appendix}
\begin{center}
{{\bf Table I a: } Lie algebras $\mf k$ s.t. $\R+\mf k$ gives rise to irreducible multiplicity free actions}
\end{center}
$$
\begin{array}{|lrclr|}
\hline
{\mf{su}(n)} & n \geq 1 &\quad\quad & {\mf{so}(n)} & n \geq 3 \\
{\mf{sp}(n)} & n \geq 2 & & {S^2(\mf{su}(n)) } & n \geq 2 \\
{\Lambda^2 (\mf{su}(n)) } & n \geq 4 & & {\mf{su}(n) \otimes
\mf{su}(m)} &
n,m \geq 2 \\
{\mf{su}(2) \otimes \mf{sp}(n)} $\ \ \ \ $ & n \geq 2 & &

{\mf{su}(3) \otimes \mf{sp} (n)} $\ \ \ \ $&
n \geq 2 \\
{\mf{su}(n) \otimes \mf{sp}(2)} & n \geq 4
& & {\mf{spin} (7)} &  \\
{\mf{spin} (9)} &  & & {\mf{spin} (10)} &  \\
{\mf{g}_2} & n \geq 1 & & {\mf{e}_6} & n \geq 3 \\ \hline
\end{array}
$$
$ \ $ \\
\begin{center}
{{\bf Table I b:} Irreducible coisotropic actions in which the
scalars are removable}
\end{center}
$$
\begin{array}{|lrclr|} \hline
\mf{su} (n) & n \geq 2 & \quad\quad& \asp (n) & n \geq 2 \\
\Lambda^2 ( \mf {su} (n)) & n \geq 4&
 & \mf {su} (n) \otimes \mf {su} (m) & n,m \geq 2,\ n \neq m \\
\mf{spin} (10) & & & \mf {su} (n) \otimes \asp (2) & n \geq 5
\\ \hline
\end{array}
$$
$\ $ \\
\begin{center}
{{\bf Table II a:} Indecomposable coisotropic actions in which the
scalars can be removed or reduced}
\end{center}
$$
\begin{array}{|lr|}
\hline {\mf{su}(n) \mioplus[\mf{su}(n)] \mf{su}(n)}
& n \geq 3,\ a \neq b  \\
{\mf{su}(n)^* \mioplus[\mf{su}(n)] \mf{su}(n)}
& n \geq 3\ a \neq -b \\
{\mf{su}(2m) \mioplus[\mf{su}(2m)] \Lambda^2 (\mf{su}(2m)) }
& m \geq 2, \ b \neq 0  \\
{\mf{su}(2m+1) \mioplus[\mf{su}(2m+1)] \Lambda^2 (\mf{su}(2m+1)) }
& m \geq 2, \ a \neq -mb   \\
{\mf{su}(2m)^* \mioplus[\mf{su}(2m)] \Lambda^2 (\mf{su}(2m)) }
& m \geq 2, \ b \neq 0 \\
{\mf{su}(2m+1)^* \mioplus[\mf{su}(2m+1)] \Lambda^2 (\mf{su}(2m+1))
}
& m \geq 2, \ a \neq mb \\
{\mf{su}(n) \mioplus[\mf{su}(n)] (\mf{su}(n) \otimes \mf{su}(m)) }
& 2 \leq n < m,\ a \neq 0   \\
{\mf{su}(n) \mioplus[\mf{su}(n)] (\mf{su}(n) \otimes \mf{su}(m)) }
& m \geq 2,\  n \geq m+2, a \neq b    \\
{\mf{su}(n)^* \mioplus[\mf{su}(n)] (\mf{su}(n) \otimes \mf{su}(m))
}
& 2 \leq n <m, a \neq 0 \\
{\mf{su}(n)^* \mioplus[\mf{su}(n)] (\mf{su}(n) \otimes \mf{su}(m))
}
& 2 \geq m, n \geq m+2, \ a \neq b \\
{(\mf{su}(2)\otimes \mf{su}(2)) \mioplus[\mf{su}(2)] (\mf{su}(2)
\otimes
\mf{su} (n)) } & n \geq 3, \ a \neq 0 \\
{(\mf{su}(n)\otimes \mf{su}(2)) \mioplus[\mf{su}(2)] (\mf{su}(2)
\otimes \mf{sp} (m)) } & n \geq 3, m\geq 4,  b \neq 0 \\\hline
\end{array}
$$
\newpage
\begin{center}
{{\bf Table II b:} Indecomposable coisotropic actions in which the
scalars cannot be removed or reduced}
\end{center}
$$
\begin{array}{|lr|}  \hline
{\mf{su}(2) \mioplus[\mf{su}(2)] \mf{su}(2)} &\\
{\mf {su}(n)^{(*)} \mioplus[\mf{su} (n)^*] (\mf{su} (n) \oplus
\mf{su} (n)) }
& n \geq 2 \\
{(\mf{su}(n+1)^{(*)} \mioplus[\mf{su}(n+1)] (\mf{su}(n+1) \otimes
\mf{su}(n)) }
& n \geq 2  \\
{(\mf{su}(2) \mioplus[\mf{su}(2)] (\mf{su}(2) \otimes \mf{sp}(m))
}
& m \geq 2 \\
{(\mf{su}(2) \oplus \mf{su}(2)) \mioplus[\mf{su}(2)] (\mf{su}(2) \otimes \mf{sp} (m)) } &\\
{(\mf{sp} (n) \oplus \mf{su}(2)) \mioplus[\mf{su}(2)] (\mf{su}(2)
\otimes \mf{sp} (m)) }
& n,m \geq 2 \\
{\mf{sp} (n) \mioplus[\mf{sp} (n)] \mf{sp} (n)}
& n \geq 2 \\
{\mf{spin} (8) \mioplus[\mf{sp} (8)] \mf{so}(8)} & \\ \hline
\end{array}
$$
\begin{small}
\n In the previous Tables we use the notation of \cite{BR}, as an
example $\su(n)\oplus_{\su(n)}{\su(n)}$  denotes the Lie
algebra $\su(n)$ acting on $\C^n \oplus \C^n$ via the direct sum
of two copies of the natural representation.
\end{small}
$\ $ \\
%massimali
\begin{center}
{{\bf Table III:} Maximal subgroups of $\SU(n)$}
\end{center}
$$
\begin{array}{|r|c|l|} \hline
i) & \SO(n) &  \\ \hline ii) & \Sp(m) & 2m=n  \\ \hline iii) &
\SUr{k}{n-k}  & 1 \leq k \leq n-1 \\ \hline iv)  & \SU(p) \otimes
\SU(q) & pq=n,\ p \geq 3,\ q \geq 3  \\ \hline v)   & \rho(H) & H
\ \rm{simple} \ \rho \in \Irr_{\C},\ deg \rho=n  \\ \hline
\end{array}
$$
$\ $ \\
\begin{center}
{{\bf Table IV:} Maximal subgroups of $\Sp(n)$}
\end{center}
$$
\begin{array}{|r|c|l|} \hline
i) & \U(n) &  \\ \hline ii) & \Sp(k) \times Sp(n-k)& 1 \geq k \geq
n-1  \\ \hline iii) & \SUr{k}{n-k}              & 1 \leq k \leq
n-1 \\ \hline iv)  & \SO(p) \otimes \Sp(q) & pq=n,\ p \geq 3,\ q
\geq 1  \\ \hline v)   & \rho(H) & H \ \rm{simple} \ \rho \in
\Irr_{\H},\ deg \rho=2n  \\ \hline
\end{array}
$$
\newpage
\begin{center}
{ {\bf Table V:} Reduced prehomogeneous triplets}
\end{center}
$$
\begin{array}{|lllr|} \hline
\quad G&\quad \rho&\quad V&\text{conditions}\\\hline
\GL (n) &  \ 2 \Lambda_1 &  \ V( \frac{n(n+1)}{2} ) & n \geq 2 \\
\GL (2n) &  \Lambda_2 &  \ V( n(n-3)) & n \geq 3 \\
\GL (2) &  3 \Lambda_1 &  \ V(4)  &\\
\GL (6) &  \Lambda_3 &  \ V(20)  &\\
\GL (7) &  \Lambda_3 &  \ V( 35 ) &\\
\GL (8) &  \Lambda_3 &  \ V(56)  &\\
\SL (3) \times \GL (2) & \ 2 \Lambda_1 \otimes \Lambda_1 &
\ V(6) \otimes V(2)  &\\
\SL (6) \times \GL (2) & \  \Lambda_2 \otimes \Lambda_1 &
\ V(15) \otimes V(2)  &\\
\SL (5) \times \GL (3) & \ \Lambda_2 \otimes \Lambda_1 &
\ V(15) \otimes V(3)  &\\
\SL (5) \times \GL (4) & \  \Lambda_2 \otimes \Lambda_1 &
\ V(10) \otimes V(4)  &\\
\SL (3) \times \SL (3) \times \GL (2) & \ \Lambda_1 \otimes
\Lambda_1 \otimes \Lambda_1 &
\ V(3) \otimes V(3) \otimes V(2)  &\\
\Sp (n) \times \GL (2m) & \ \Lambda_1 \otimes \Lambda_1 &
\ V(2n) \otimes V(2m) & n \geq 2m \geq 2 \\
\GL (1) \times \Sp (3) & \  \rho_1 \otimes \Lambda_3 &
\ V(1) \otimes V(14)  &\\
\SO (n) \times \GL (m) &  \  \Lambda_1 \otimes \Lambda_1 &
\ V(n) \otimes V(m) & n \geq 3 , \frac{n}{2} \geq m \geq 1 \\
\GL (1) \times \Spin  (7) &  \ \rho_1 \otimes {\rm spin \ rep.} &
\ V(1) \otimes V(8)  &\\
\Spin (7) \times \GL (2) &  \ {\rm spin \ rep} \otimes \Lambda_1 &
\ V(8) \otimes V(2)  &\\
\Spin (7) \times \GL (3) &  \ {\rm spin \ rep} \otimes \Lambda_1 &
\ V(8) \otimes V(3)  &\\
\GL (1) \times \Spin (9) & \ \rho_1\otimes {\rm spin \ rep.} &
\ V(1) \otimes V(16)  &\\
\Spin (10) \times \GL (2) & \ {\rm half-spin \ rep.} \otimes
\Lambda_1 &
\ V(16) \otimes V(2)  &\\
\Spin (10) \times \GL (3) & \ {\rm half-spin \ rep.} \otimes
\Lambda_1 &
\ V(16) \otimes V(3)  &\\
\GL (1) \times \Spin (11) & \ \rho_1 \otimes {\rm spin. \ rep.} &
\ V(1) \otimes V(32)  &\\
\GL (1) \times \Spin (11) & \ \rho_1 \otimes {\rm half-spin \ rep.} &
\ V(1) \otimes V(64)  &\\
\GL (1) \times \G_2  & \ \rho_1 \otimes \Lambda_2 &
\ V(1) \otimes V(7)  &\\
\G_2 \times \GL (2) & \  \Lambda_2 \otimes \Lambda_1 &
\ V(7) \otimes V(2)  &\\
\GL (1) \times \Ea_6  & \ \rho_1 \otimes \Lambda_1 &
\ V(1) \otimes V(27)  &\\
\Ea_6  \times \GL (2) & \  \Lambda_1 \otimes \Lambda_1 &
\ V(27) \otimes V(2)  &\\
\GL (1) \times \Ea_7  & \  \rho_1\otimes \Lambda_1 &
\ V(1) \otimes V(56)  &\\
\Sp (n)   \times \GL (2) & \  \Lambda_1 \otimes 2 \Lambda_1 &
\ V(2n) \otimes V(3)  &\\
\GL (1) \times \Sp (n) \times \SO (3) & \  \rho_1 \otimes \Lambda_1 \otimes
\Lambda_1 & \ V(1) \otimes V(2n) \otimes V(3)  &\\
\SL (n)  \times \GL (m) & \  \Lambda_1 \otimes \Lambda_1 &
\ V(n) \otimes V(m) & \frac{m}{2} \geq n \geq 1  \\
\GL (2m +1) & \  \Lambda_2 &  \ V(m(2m+1)) & m \geq 2  \\
\SL (2m+1)  \times \GL (2) & \  \Lambda_2 \otimes \Lambda_1 &
\ V(m(2m+1)) \otimes V(2) &   m \geq 2\\
\Sp (n)   \times \GL (2m+1)) & \  \Lambda_1 \otimes \Lambda_1 &
\ V(2n) \otimes V(2m+1) & n>2m+1 \geq 1 \\
\GL (1)  \times \Spin (10) & \  \rho_1 \otimes {\rm half-spin \ rep.} &
\ V(1) \otimes V(16)  &\\ \hline
\end{array}
$$
\begin{small}
\n In the previous Table  $\rho_1$ and $\Lambda_i$  denote the standard representation of $\GL(1)$ and  the fundamental highest weights respectively. 
\end{small}

$\ $ \\
\n
Leonardo Biliotti \\
Dipartimento di Matematica e Applicazioni per l'Architettura, \\
Universit\`a di Firenze, Piazza Ghiberti 27, 50142 Firenze, Italy. \\
email: biliotti@math.unifi.it \\
$\ $ \\
Anna Gori \\
Dipartimento di Matematica ``Ulisse Dini'', \\
Universit\`a di Firenze, Viale
Morgagni 67-A, 50134 Firenze, Italy. \\
email: gori@math.unifi.it

\begin{thebibliography}{99}
\bibitem{BR}
{\sc Benson, C., Ratcliff, G.:}
\newblock {\em A classification of multiplicity free actions},
\newblock  J. Algebra,
{\bf 181}, 152--186 (1996).
%
\bibitem{Bergmann}
{\sc Bergmann I.:}
\newblock {\em Reducible polar representations,}
\newblock  Manus. Math.,
{\bf 104}, 309--324 (2001).
%
\bibitem{Da}
{\sc Dadok J.:}
\newblock {\em Polar cordinates induced by actions of compact Lie Groups,}
\newblock  Trans. Am. Math. Soc.,
{\bf 288}, 125--137 (1985).
%
%
\bibitem{GS}
{\sc Guillemin, V. Sternberg, S.:}
\newblock {\em Symplectic Techniques in Physics}
\newblock Cambridge: Cambridge University Press, 1984.
%
\bibitem{He}
{\sc Heintze, E. Eschemburg J.:}
\newblock {\em On the classification of Polar representation,}
\newblock  Math. Z.,
{\bf 286}, 391--398 (1999).
%
\bibitem{Hel}
{\sc Helgason, S.:}
\newblock {\em Differential Geometry, Lie Groups and Symmetric Spaces}
\newblock  New York-London: Academic Press-inc, 1978.
%
%
\bibitem{Hw}
{\sc Huckleberry, A.T., Wurzbacher, T.:}
\newblock {\em Multiplicity-free complex manifolds,}
\newblock  Math. An\-nalen.,
{\bf 286}, 261--280 (1990).
%
\bibitem{Kac}
{\sc Kac, V.G.:}
\newblock {\em Some remarks on Nilpotent orbit,}
\newblock J. Algebra,
{\bf 64}, 190--213 (1980).
%
\bibitem{Kol} {\sc Kollross, A.:}
\newblock {\em A classification of hyperpolar and cohomogeneity one actions,}
\newblock Trans. Am. Math. Soc.,
{\bf 354}, 571--612 (2002).
%
\bibitem{Kr}
{\sc Kraft, H.:}
\newblock {\em Geometrische Methoden in der Invariantentheorie,}
\newblock Braunschweig, Wiesbaden: Vieweg 1984.
%
\bibitem{On}
{\sc Onishchik, A. L.:}
\newblock {\em Inclusion relations among transitive compact transformation
groups,}
\newblock Trans. Am. Math. Soc.,
{\bf 50}, 5--58 (1966).
%
%
\bibitem{PoT}
{\sc Podest\`a, F., Thorbergsson, G.:}
\newblock {\em Polar and {C}oisotropic {A}ctions on {K}\"ahler
{M}anifolds,}
\newblock Trans. Am. Math. Soc.,
{\bf 354}, 236--238 (2002).
%
%
\bibitem{SK}
{\sc Sato, M., Kimura, T.:}
\newblock {\em A classification of irreducible prehomogeneous vector spaces
and their relative invariants,}
\newblock Nagoya Math.  J.,
{\bf 65}, 1--155 (1977).
%
%\begin{thebibliography}{99}
%\addcontentsline{toc}{chapter}{References}
%\pagestyle{myheadings}
%\markright{References}
\end{thebibliography}
\end{document}